\documentclass[onecolumn,11pt]{article}
\usepackage{algorithmic}
\usepackage[ruled,vlined]{algorithm2e}
\usepackage{xcolor}
\usepackage{caption}
\usepackage[utf8]{inputenc}
\usepackage{amsmath}

\usepackage[colorinlistoftodos]{todonotes}

\usepackage{soul}
\definecolor{dr}{rgb}{0.75,0.00,0.00}
\definecolor{lr}{rgb}{1.00,0.75,0.75}
\sethlcolor{lr}

\usepackage{multirow}
\usepackage{graphicx}
\usepackage{amsmath,amsfonts}
\usepackage{breqn}
\usepackage{subcaption}
\usepackage{float}
\usepackage{enumerate}
\usepackage{booktabs}
\restylefloat{table}
\usepackage[a4paper,margin=1in]{geometry}
\usepackage[onehalfspacing]{setspace}
\usepackage[utf8]{inputenc}
\usepackage{enumitem}
\usepackage{verbatim}
\usepackage[font=small,labelfont=bf]{caption}
\usepackage{tabularx}
\usepackage{makecell}
\usepackage{bbm}
\usepackage{url}
\usepackage{enumerate}
\usepackage{authblk}
\usepackage[short,nocomma]{optidef_2020} 
\usepackage[authoryear,round]{natbib}

\usepackage[colorlinks = true,
linkcolor = blue,
urlcolor = blue,
citecolor = black,
anchorcolor = blue]{hyperref}

\DeclareMathOperator*{\lexmin}{lex~min}

\title{The Impact of Congestion and Dedicated Lanes \\ on On-Demand Multimodal Transit Systems}
\author{Jason~Lu}
\author{Anthony~Trasatti}
\author{Hongzhao~Guan}
\author{Kevin~Dalmeijer}
\author{Pascal~Van~Hentenryck}
\affil{H. Milton Stewart School of Industrial and Systems Engineering,\protect\\ Georgia Institute of Technology}
\date{\today}

\begin{document}

\maketitle


\vfill

\begin{abstract}
\noindent
Traffic congestion can have a detrimental effect on public transit systems, and understanding and mitigating these effects is of critical importance for effective public transportation.
Implementing Dedicated Bus Lanes (DBLs) is a well-known intervention to achieve this goal.
A DBL is a designated
lane for bus transit, which avoids congestion and substantially lowers
the travel time.  This makes transit more attractive,
encouraging more travelers to adopt public transportation.
This paper studies the impact of congestion and DBLs on novel On-Demand Multimodal Transit Systems (ODMTS).
ODMTS combine
traditional rail and bus networks with on-demand
shuttles.  Previous case studies have shown that ODMTS may
simultaneously improve travel time, reduce system cost, and attract
new passengers.  Those
benefits were shown for an ideal world without traffic congestion, and
this paper hypothesizes that the advantages of ODMTS can be even more
pronounced in the real world.  This paper explores this hypothesis by
creating realistic congestion scenarios and solving bilevel
optimization problems to design ODMTS under these scenarios.  The
impact of DBLs on ODMTS is evaluated with a comprehensive case
study in the Metro Atlanta Area. The results show that DBLs can
significantly improve travel times and are effective at increasing
adoption of the system.\\

\noindent\emph{\textbf{Keywords}:
Congestion, Dedicated Bus Lanes, Public Transit, On-Demand Multimodal Transit Systems.
}\\
\end{abstract}

\clearpage
\section{Introduction}
Traffic congestion can have a detrimental effect on public transit systems, and understanding and mitigating these effects is of critical importance for effective public transportation.
Implementing Dedicated Bus Lanes (DBLs) is a well-known intervention to achieve this goal.
A DBL is a designated lane for bus transit, which avoids congestion and substantially lowers the travel time for riders. This makes transit more attractive, encouraging more travelers to adopt public transportation.
Increased ridership in turn leads to a plethora of benefits, including fewer cars on the road, less emission, and increased revenue for transit operators that can be used to further improve service.

\begin{figure}[!t]
	\centering
	\includegraphics[width=0.7\linewidth]{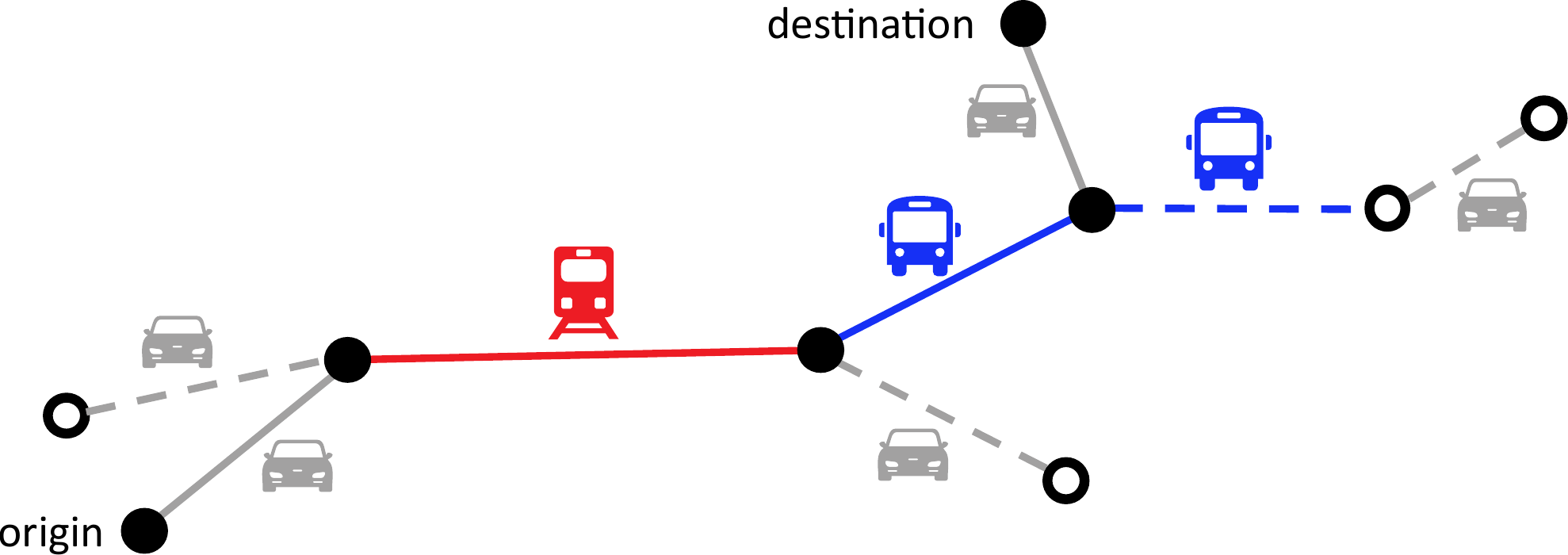}
	\caption{Example ODMTS with Passenger Path (solid lines)}
	\label{fig:odmts_example}
\end{figure}

This paper studies the impact of congestion and DBLs on novel On-Demand Multimodal Transit Systems (ODMTS).
ODMTS combine traditional rail and bus networks with
on-demand shuttles \citep{VanHentenryck2019-DemandMobilitySystems}.
Trains and buses serve the busy corridors on a fixed schedule, while
shuttles dynamically serve the first and last miles.
Figure~\ref{fig:odmts_example} provides an example of an ODMTS and the
path of a single passenger.  Passengers provide their origin and
destination through a mobile application, after which the ODMTS
provides a route to serve them.  In this case, a passenger is picked
up by an on-demand shuttle close to their origin and brought to the
train station.  The passenger is then instructed to take a train and a
bus, which both run on a fixed schedule.  When the passenger arrives
at the bus station closest to the destination, another on-demand
shuttle is ready to pick them up and serve the last mile.  The
\citet{SAML2020-HowDemandMultimodal} provides a
\href{https://sam.isye.gatech.edu/projects/demand-multimodal-transit-systems}{video}
of this process.

Case studies in Canberra, Australia \citep{maheo2019benders}, Ann
Arbor, Michigan
\citep{basciftci2020bilevel,BasciftciVanHentenryck2022-CapturingTravelMode,AuadVanHentenryck2022-RidesharingFleetSizing},
and Atlanta, Georgia (\citet{dalmeijer2020transfer,auad2021resiliency,GuanEtAl2022-HeuristicAlgorithmsIntegrating}; \citet{Vanhentenryck2023reach})
demonstrate that ODMTS may simultaneously improve travel time, reduce
system cost, and attract new passengers compared to the existing
systems.  Those benefits are shown for an ideal world without traffic
congestion, and this paper hypothesizes that the advantages of ODMTS
can be even more pronounced in the real world.  By the very nature of
ODMTS, shuttle trips are inherently local and minimally affected by
traffic.
Furthermore, buses are only used to serve high-density corridors, making DBLs for ODMTS potentially more impactful than for a traditional system.

\begin{figure}[t]
	\centering
	\includegraphics[width=0.7\linewidth]{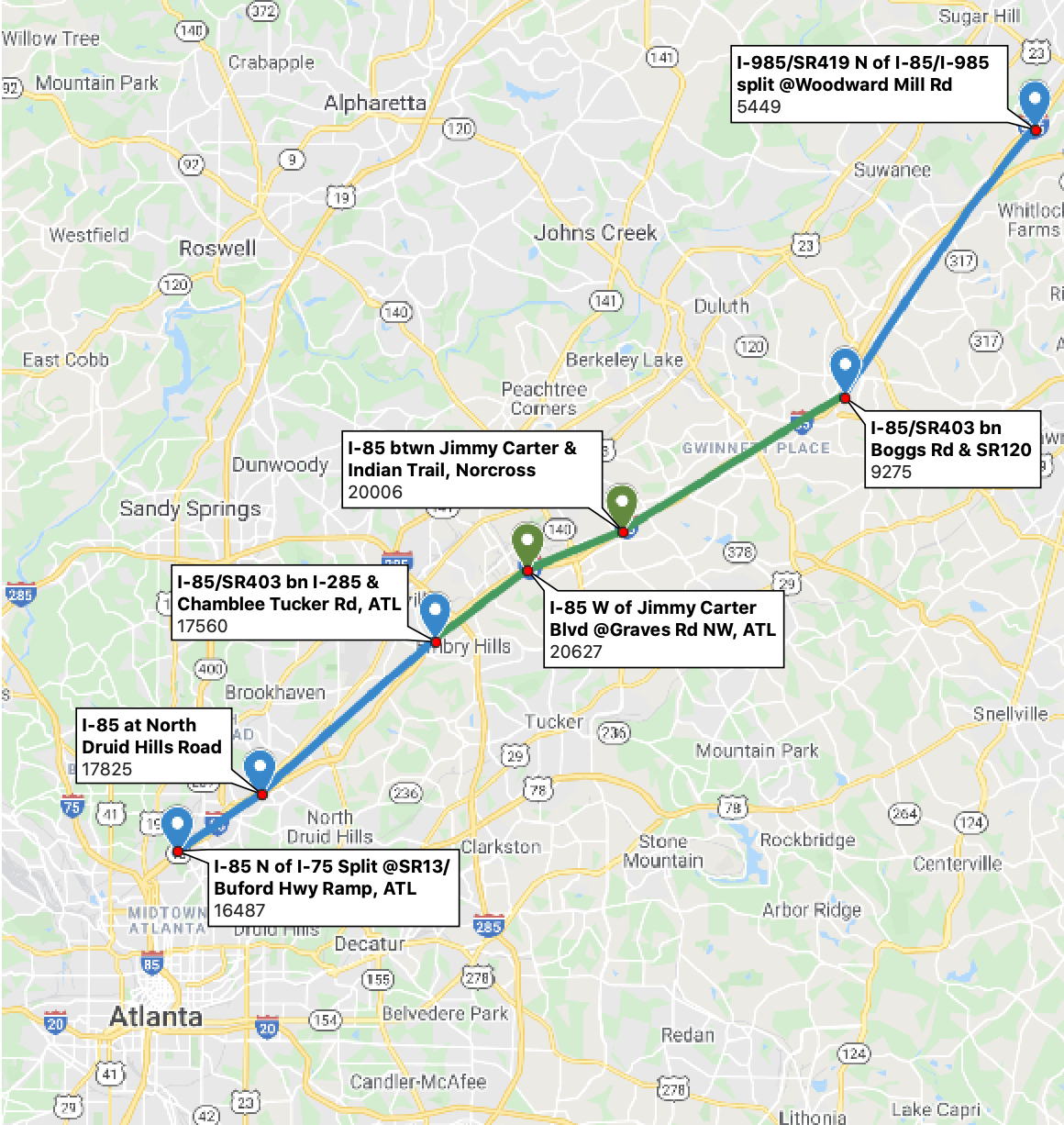}
	\caption{Average Traffic Volume along I-85 during the 7am-9am Rush Hour (March 21-25, 2022)}
	\label{fig:gdot_sensors}
\end{figure}

A case study is conducted to fill the critical gap in understanding the impact of congestion on ODMTS.
Furthermore, it is shown that DBLs may lead to significantly increased rider adoption.
The
case study focuses on traffic from Gwinnett Country to the city of
Atlanta, both in the Metro Atlanta Area (Georgia, USA).  Gwinnett
County is the second-most populated county in Georgia, with close to a
million residents \citep{CensusBureau2021-AnnualEstimatesResident}.
Many of these residents work in Atlanta and drive the Interstate 85
(I-85) for their commute, creating significant traffic.
Figure~\ref{fig:gdot_sensors} summarizes Average Traffic Volume (ATV)
from continuous count stations along I-85
\citep{GDT2022-TrafficAnalysisData}.  When driving south from Gwinnett
to Atlanta, the ATV more than doubles on the highlighted section,
which shows the importance of this highway and the high potential for
congestion.  Table~\ref{tab:google_maps_info} uses data from Google
Maps to demonstrates the enormous potential of DBLs to reduce travel
time on I-85.  The table compares rush hour travel time to free flow
travel time for both the full segment and for the highlighted segment.
It can be seen that a DBL may save commuters from Gwinnett up to 61
minutes every morning.

\begin{table}[t]
	\centering
	\begin{tabular}{lll}
		\toprule
		& Full segment & Highlighted segment\\
		\midrule
		Distance & 28.8 miles & 12.0 miles \\
		Rush hour time (7am-9am) & 28-85 min & 12-50 min \\
		Free flow time (2am) & 24 min & 10 min \\
		Potential savings & 4-61 min & 2-40 min \\
		\bottomrule
	\end{tabular}%
	\caption{Travel Times and Potential Savings on I-85 (March 21-25, 2022)}
	\label{tab:google_maps_info}
\end{table}

The impact of DBLs is studied by creating different congestion
scenarios and designing ODMTS with and without DBLs.  Designing the
optimal ODMTS is modeled as a bilevel optimization problem in which
passengers choose whether to adopt the system based on the quality of
the trip they are offered.  This work combines the general model by
\citet{auad2021resiliency} with the adoption model by
\citet{BasciftciVanHentenryck2022-CapturingTravelMode}.  To solve this
model at the scale of the case study, this paper applies the heuristic
algorithm by \citet{GuanEtAl2022-HeuristicAlgorithmsIntegrating} and
uses Benders decomposition to solve the fixed-demand ODMTS design
problem in every iteration.  This enables a comprehensive case study
in the Metro Atlanta Area that considers different levels of
congestion, shows the impact of DBLs, and demonstrates the increase in
adoption of the ODMTS.

The contributions of the paper can be summarized as follows:

\begin{enumerate}
\item The paper creates realistic congestion scenarios and solves bilevel optimization problems to design ODMTS under these scenarios.
\item The paper evaluates the effect of congestion on ODMTS with a comprehensive case study in the Metro Atlanta Area.
\item The paper demonstrates on the case study that DBLs can significantly increase ODMTS adoption.
\end{enumerate}

The remainder of this paper is structured as follows.
Section~\ref{sec:literature-rev} gives an overview of the literature.
Section~\ref{sec:methodology} presents the methods for generating
congestion scenarios and designing ODMTS with and without DBLs.
Section~\ref{sec:current} introduces the case study area and the
experimental settings.  The effect of congestion on ODMTS without DBLs
is analyzed in Section \ref{sec:baseline}, after which Section
\ref{sec:impact-of-DBLs} analyzes the effect of adding DBLs to the
system.  Section \ref{sec:discussion} concludes the study with final
thoughts on DBLs in ODMTS.

\section{Literature Review}
\label{sec:literature-rev}
DBLs have been studied through a variety of simulations and case
studies.  \citet{BassoEtAl2011-CongestionPricingTransit} analyse a
model where travelers have the choice between car, bus, and an outside
option.  The authors optimize frequency, vehicle size, spacing between
stops, and capacity on DBLs to study different urban congestion
management policies.  It is found that implementing DBLs is a better
stand-alone policy than subsidizing transit or pricing congestion.
\citet{RussoEtAl2022-DedicatedBusLanes} study the welfare effect of
DBLs in Rome, Italy and conclude that DBLs may decrease bus travel
time by 18\%, decrease waiting time by 12\%, and increase ridership by
26\%.  Furthermore, travel time for other motor vehicles may improve
as well.  Other case studies include the study by
\citet{Ben-DorEtAl2018-AssessingImpactsDedicated} for Sioux Falls,
South Dakota.  The authors present an agent-based model and use the
MATSim simulator to demonstrate that DBLs fundamentally improve the
effectiveness of public transit and essentially make public transit trip durations similar during peak hours and non-peak hours.
\citet{HoonsiriEtAl2021-UsingCombinedBus} consider an intervention
where fixed buses are allowed onto dedicated Bus Rapid Transit (BRT)
lanes.  A case study on the Rama 3 road in Bangkok, Thailand shows
potential savings in travel time and a reduction in greenhouse gas
emissions. \citet{Tsitsokas2021modeling} formulate a nonlinear combinatorial optimization model to solve a DBL allocation problem on a large-scale road network. A case study in the San Francisco Metropolitan area demonstrates that the road network with DBL configurations significantly improves travel times for both cars and bus users compared to no DBLs.

Dedicated lanes have also been studied in the context of Autonomous
Vehicles (AVs), where similar benefits are observed.  \citet{ChakrabortyEtAl2021-FreewayNetworkDesign}
consider the problem of deciding which lanes on a freeway network
should be AV-exclusive.  Traveler choice between AVs and regular
vehicles is captured by a logit model, and the resulting non-convex
mixed-integer nonlinear program is solved through a combination of
Benders decomposition and the Method of Successive Averages.  The
authors construct multiple freeway networks with up to 23 links to
demonstrate the algorithm.
\citet{ChenEtAl2020-ModelingControlAutomated} explore the idea of
giving AVs access to dedicated BRT lanes and focus on the interaction
between AVs and BRTs.  The authors first consider a model for
mixed-use lanes, and propose a sequential optimization method to
analyze the performance.  The second model allows AVs to move in and
out of the BRT lane at different locations, and is solved through
successive linear programming.  A case study for BRT Line 1 in
Beijing, China shows that mixed-use lanes can both improve efficiency
of AVs and reduce congestion on the other lanes.

The ODMTS concept was proposed by \citet{maheo2019benders}. They leverage ideas from 
  hub arc location problems for transportation and telecommunication networks \citep{CampbellEtAl2005-HubArcLocation} to design a hub and shuttle transit system that determines which bus arcs to open and how to serve the first/last miles with shuttles.
A Benders decomposition algorithm is introduced to solve the problem, and a case study in Canberra, Australia demonstrates that the
new system improves transit times without negatively affecting system costs.
\citet{auad2021resiliency} build on this work to present an end-to-end ODMTS solution that combines demand estimation, network design, fleet sizing, and real-time shuttle dispatching.
Earlier methods are adapted to the multimodal setting, and practical constraints are introduced, such as a limit on the numbers of passenger transfers.
The pipeline is used to conduct a case study in Atlanta, Georgia and to demonstrate the resiliency of ODMTS under various scenarios of COVID-19 pandemic response.

Recent studies on ODMTS have focused on including latent demand into the design of the network.
\citet{basciftci2020bilevel, BasciftciVanHentenryck2022-CapturingTravelMode} study ODMTS Design with Adoptions (ODMTS-DA) with bilevel optimization models.
The ODMTS-DA aims to design an ODMTS while taking into account that new passengers may adapt the system based on a personalized choice model.
Exact methods are introduced to solve the ODMTS-DA, and a small-scale case study in Ann Arbor, Michigan demonstrates that ODMTS can improve service for the existing riders and attract new riders at the same time.
\citet{GuanEtAl2022-HeuristicAlgorithmsIntegrating} develop heuristics to solve the ODMTS-DA for large-scale systems to which the exact method does not scale.
These heuristics solve the problem iteratively by solving an ODMTS design problem with fixed demand, evaluating the passenger choices based on the new network, updating the ODMTS design problem, and repeating this process until convergence.
A case study in Atlanta, Georgia is used to demonstrate that the heuristic is able to find high-quality solutions in significantly less time than the exact algorithm.

Other studies have investigated the coordination between on-demand services and public transit systems.
\citet{SalazarEtAl2018-InteractionAutonomousMobility} develop a network flow model and a pricing scheme to capture the interaction of Autonomous Mobility-on-Demand (AMoD) fleets with public transit.
They undertake a case study on the existing transit system in New York City and show that integrating AMoD fleets with public transit brings significant benefits compared to AMoD fleets operating independently.
\citet{PintoEtAl2020-JointDesignMultimodal} study the allocation of resources between transit systems and shared-use autonomous mobility services.
They design a bilevel mathematical programming formulation that decides the balance between open bus routes and fleet sizes.
Their case study in Chicago demonstrates significant benefits to passenger waiting times compared to the existing multimodal transit network.

The current paper studies whether the advantages of DBLs translate to
ODMTS.
Previous works on DBLs either study the effect of policy
changes on existing systems or optimize how to modify existing
systems.  This paper leverage the fact that ODMTS is a novel type of
transit system to introduce models that take congestion into account
\emph{while designing the system}.  It also presents a case study for
a large multi-modal system, while related optimization-based studies
often focus on a single mode and a single corridor.

\section{Methodology}
\label{sec:methodology}

This section introduces the methods that are needed to rigorously
study the effect of congestion on ODMTS and to assess the benefit of
DBLs.  It first proposes an approach to model congestion in ODMTS with limited resources, 
and how to create the congestion scenarios that are used in the analysis. 
It then presents a bilevel optimization model to design ODMTS for these
scenarios. It concludes by discussing how to incorporate DBLs into
the design.

\subsection{Congestion Modeling}
\label{subsec:congestion-modeling}

Obtaining travel time predictions for various levels of congestion for
each rider becomes computationally challenging or financially
expensive for a large ridership.  Therefore, this paper constructs
congestion scenarios in two steps.  The first step is to obtain a
basis of travel times, and the second step is to multiply the travel
time basis with a factor that depends on the congestion scenario and
the local conditions.

\paragraph{Step 1: Travel Time Basis}
\label{sssec:o-d-travel-time-basis}

The travel time basis serves as the default travel time for any origin-destination pair (OD pair) in this study. For a given congestion scenario, the travel time basis is multiplied by a congestion scaling factor, resulting in the travel times under that scenario.
The travel time basis is obtained from the POLARIS Transportation
System Simulation Tool \citep{AuldEtAl2016-PolarisAgentBased}. 
POLARIS is an agent-based modeling software tool created by Argonne
National Laboratory that uses a dynamic simulation of travel demand,
network supply, and network operations.  It is used for a variety of
cities to simulate an average day of activities, and has been tuned
specifically for Atlanta using survey data from the Atlanta Regional
Commission \citep{ARC2022-ActivityBasedModeling}.  POLARIS provides
travel times between every pair of Traffic Analysis Zones (TAZs).  The
travel time between arbitrary points is approximated by the travel
time between the corresponding TAZs. POLARIS \emph{does} model
congestion, but tuning to the survey data only provides a single scenario.
This motivates the next step.

\paragraph{Step 2: Congestion Scaling Factors}
\label{sssec:congestion-variation}

\begin{figure}[t]
	\centering
	\includegraphics[width=.7\linewidth, height=.7\linewidth]{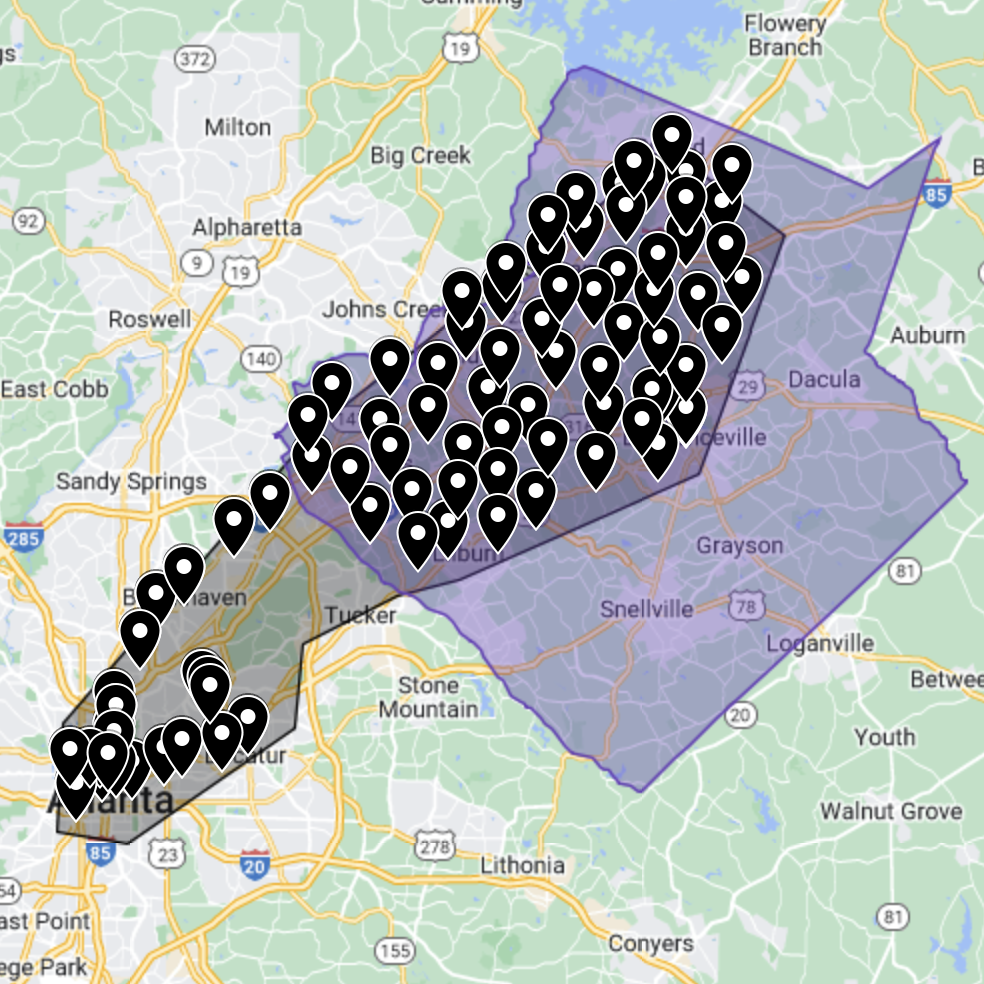}
	\caption{Query Reference Locations with Respect to the Case Study Area (Gray) and Gwinnett County (Purple)}
	\label{fig:hubs_set1}
\end{figure}

To determine an appropriate congestion scaling factor for every OD pair (without querying all of them), this paper introduces a set $Q$
of Query Reference Locations (QRL).  The set $Q$ is shown by
Figure~\ref{fig:hubs_set1}, and consists of existing Gwinnett bus
stops, MARTA rail stations, and local points of interest in Gwinnett.
For a given OD pair and a scenario $k$, the idea is to first identify the reference location
$i\in Q$ closest to the origin and $j \in Q$ closest to the
destination.  The congestion scaling factor $R_{ij}^k$ is then
calculated as $\tau_{ij}^k/\tau_{ij}^0$, where the travel time
$\tau_{ij}^k$ for scenario $k$ is queried and $\tau_{ij}^0$ is
from the travel time basis in Step 1.  The estimated travel time for the OD pair under scenario $k$ is obtained by
multiplying the OD travel time basis by $R_{ij}^k$.  The benefit of
this approach is that it allows for region-specific scaling but only
queries travel times between the QRLs, which avoids prohibitively
expensive data collection.  In the case where $i = j$, i.e., the
origin and the destination are mapped to the same QRL, the scaling factor
is calculated between $i$ and the QRL closest to $i$.

\paragraph{Scenarios}

The baseline scenario $k=0$ is provided by the POLARIS travel time
basis.  To generate the congestion scenarios, this paper uses the
Directions API on \citet{directionsAPI} to obtain
scaling factors for different levels of congestion.  The Directions
API supports two traffic models: \texttt{best\_guess} and
\texttt{pessimistic}, which approximate travel times under congestion.  Based on these traffic models, three congestion
scenarios are created:
\begin{itemize}[noitemsep]
	\item Expected: use \texttt{best\_guess} travel times.
	\item 50-50: use the average of \texttt{best\_guess} and \texttt{pessimistic} travel times.
	\item Pessimistic: use \texttt{pessimistic} travel times.
\end{itemize}
The corresponding travel times $\tau_{ij}^{exp}$, $\tau_{ij}^{50}$, and $\tau_{ij}^{pes}$ are queried between all $i,j \in Q.$ Then the scaling factors $R_{ij}^{exp}$, $R_{ij}^{50}$, and $R_{ij}^{pes}$ between all $i,j \in Q$ are calculated accordingly and define the scenarios.

\subsection{ODMTS Design}
\label{subsec:odmts-design}

This paper combines the general model by \citet{auad2021resiliency}
with the adoption model by
\citet{BasciftciVanHentenryck2022-CapturingTravelMode} to design an
ODMTS that takes passenger adoption into account.  The resulting
bilevel model inherits support for transfer constraints, multiple
frequencies, and both bus and rail connections.
A full description is included below to clarify the assumptions and make the paper more self-contained.
A more implicit description of a similar model can be found in \citet{GuanEtAl2022-HeuristicAlgorithmsIntegrating}.

\paragraph{Problem Description}
In the general model from \citet{auad2021resiliency}, the first input to the ODMTS network design problem is a transit network modeled by the directed multigraph $G=(V,A)$.
Vertices $V$ represent locations, and arcs $a \in A$ capture potential methods of transportation between these locations.
Every arc from $i(a)\in V$ to $j(a)\in V$ is associated with a mode $m(a) \in \{shuttle, bus, rail\}$ that indicates how this connection is served, and bus and rail arcs additionally have a frequency $f(a)$
from the set $F^{bus}$ or $F^{rail}$ respectively.
Bus and rail arcs are used to transport riders between a set of transit hubs $V_H \subseteq V$, while shuttle arcs provide the connections to and from the hub network.
For example, arc $a \in A$ may correspond to opening a bus connection with frequency 10 per hour between hub $i \in V_H$ and hub $j \in V_H$.
It is important to note that $G$ may have parallel arcs that model the same connection, but differ in mode or frequency.
This is the case when a bus service can be offered at different frequencies, for example.
The network design problem decides which potential arcs $a \in A$ to make available.
In line with earlier work it is assumed that shuttle connections and the fixed rail network are always available, and the model only decides whether or not to open each bus arc \citep{auad2021resiliency,GuanEtAl2022-HeuristicAlgorithmsIntegrating}.
Furthermore, bus and rail capacities are ignored, as the high-frequency vehicles are almost never full in practice \citep{auad2021resiliency}.
Shuttle capacities are captured implicitly: the cost of using a shuttle represents the average cost for serving passengers in a ridesharing system with small vehicles.
For convenience, the set of bus arcs is denoted by $A^{bus}$.
Finally, every arc $a \in A$ is associated with a distance $d_a$, a travel time $\tau_a$, and an expected waiting time $\omega_a$.
For shuttles and buses the travel time is based on the congestion scenario (as discussed in Section~\ref{subsec:congestion-modeling}), while for rail the travel time approximates the current schedule.

The second input to the network design problem is a set $T$ of passenger trips.
Each trip $t\in T$ is an OD pair with origin $o(t)\in V$ and destination $d(t)\in V$, and is associated with $p(t)$ number of travelers.
Capturing latent demand is essential when designing public transit services, as omitting these potential riders may cause inefficiencies and unfairness \citep{BasciftciVanHentenryck2022-CapturingTravelMode}.
For this reason, the set $T$ consists of both \emph{potential riders} (i.e., latent demand) $T' \subseteq T$ who will only adopt the ODMTS if their travel time is sufficiently short, and \emph{existing riders} $T\setminus T'$ who will always adopt the ODMTS. This study assumes existing riders have no alternative travel options, because they are already transit users.
This behavior is captured by a choice function $\mathcal{C}^t$ that will be defined later.
A trip that is served by the ODMTS is modeled as a path in $G$ from $o(t)$ to $d(t)$.
To prevent a large number of transfers, these paths are limited to at most $K \ge 1 $ arcs. 

The objective is to minimize a convex combination $\alpha \in [0,1]$ of system cost (weight $1-\alpha$) and passenger inconvenience (weight $\alpha$) for the travelers who adopt the system.
The fixed cost of opening bus arc $a \in A$ is defined as $\beta_a = (1-\alpha) \tau_a f(a) c^{bus}$, which combines the travel time $\tau_a$, the number of buses over the time horizon $f(a)$, and the bus cost per hour $c^{bus}$.
The rail network is assumed to be fixed, and the corresponding constant is omitted from the objective.
Using arc $a \in A$ for trip $t \in T$ contributes to the objective as follows:
\begin{equation}
	\gamma_a^t = 
	\begin{cases}
		p(t) \left((1-\alpha) d_a \bar{c}^{shuttle} + \alpha \tau_a\right) & \textrm{ if } m(a) = shuttle\\
		p(t) \alpha \left(\tau_a + \omega_a\right) & \textrm{ if } m(a) \in \{bus, rail\}.
	\end{cases}
\end{equation}
For shuttles, the distance $d_a$ is multiplied by the shuttle cost per mile $\bar{c}^{shuttle}$, the inconvenience is the travel time $\tau_a$, and the sum is multiplied by the number of passengers $p(t)$.
For the other modes the only costs are fixed costs, such that $\gamma_a^t$ only consists of the inconvenience $\tau_a + \omega_a$, which sums the travel time and the waiting time.
For shuttles and buses, recall that $\tau_a$ is based on the congestion scenario.
This means that opening bus arcs becomes more expensive, and using shuttles and buses becomes more inconvenient as congestion increases.
Following earlier work it is assumed that shuttles are readily available and do not impose additional waiting time.
This is motivated by the fact that shuttle waiting times are typically short \citep{VanHentenryck2019-DemandMobilitySystems,auad2021resiliency} and are comparable to the time it takes to park a car or walk to a bus station, which is not penalized either.
For a given network design, passengers are always offered a route that minimizes the objective, but whether they accept this route and adopt the system depends on their personal preferences.
If trip $t\in T$ is adopted, this additionally leads to a benefit of $\zeta^t = (1-\alpha)p(t)c^{ticket}$, i.e., the number of tickets multiplied by the price per ticket.
Without loss of generality, the constant revenue from tickets sold to the existing riders is omitted from the objective.

\paragraph{Optimization Model}
The bilevel optimization model by \citet{BasciftciVanHentenryck2022-CapturingTravelMode} consists of the leader model~\eqref{formulation:leader} and a follower model~\eqref{formulation:follower} for every trip $t \in T$.
For brevity, $\delta^+(i)$ denotes the set of arcs going out of $i \in V$, and $\delta^+(i, bus)$ further restricts this set to bus arcs.
The sets $\delta^-(i)$ and $\delta^-(i,bus)$ are defined similarly for in-arcs.
The leader model uses the binary variables $z_a$ to decide which bus arcs $a\in A^{bus}$ are opened ($z_a=1$) or closed ($z_a=0$).
Constraints~\eqref{eq:masterFlowConservation} ensure that the bus frequencies are balanced at every hub, and Constraints~\eqref{eq:masterOneFrequency} ensure that at most one frequency is selected among parallel bus arcs (recall that $G$ is a multigraph and may offer multiple frequency options for the same connection).

The leader anticipates the behavior of follower $t \in T$ through a binary variable $x^t$ that represents the choice function $\mathcal{C}^t(l^t)$ (Equation~\eqref{eq:userChoiceModel}).
The choice function equals one if $t \in T$ adopts, and zero otherwise.
In particular, it is assumed that travelers adopt if they are offered a path of length $l^t$ at most equal to the adoption factor $\rho$ times the travel time by car $l^t_{car}$, i.e.,
\begin{equation}
	\mathcal{C}^t(l^t) = \mathbbm{1}\left(l^t \le \rho l^t_{car}\right).
\end{equation}
The leader Objective~\eqref{eq:objectiveOveral} then sums the fixed cost for opening bus arcs, the cost and inconvenience for the existing riders, and the cost and inconvenience including ticket revenue for the potential riders that adopt ($x^t=1$).

Based on the network design $\boldsymbol{z}$, the follower model~\eqref{formulation:follower} finds a path for trip $t \in T$ that minimizes the lexographic Objective~\eqref{eq:objectivePassenger}.
The binary variables $y_a^t$ for $a \in A$ define a unit flow from $o(t)$ to $d(t)$ that represents the path for this trip.
Constraints~\eqref{eq:extsubFlowConservation} enforce flow conservation, Constraints~\eqref{eq:extsubArcCapacity} ensure that arcs can only be used if they are opened by the leader, and Constraint~\eqref{eq:extsubTransferLimit} enforces the transfer limit.
Followers are offered a path that minimizes the primary objective of cost and inconvenience, and ties are broken by minimizing the secondary objective of trip length.
The corresponding optimal values are denoted by $g^t$ and $l^t$ respectively.
Note that travelers decide whether to adopt based on the length of the trip $l^t$, while the leaders minimize a combination of trip lengths and cost.

\begin{figure}[t]
\begin{mini!}
%
	{\boldsymbol{x}, \boldsymbol{z}}
%
	{\sum_{a \in A^{bus}} \beta_a z_a + \sum_{t \in T \setminus T'} g^t + \sum_{t \in T'} x^t \left(g^t - \zeta^t\right), \label{eq:objectiveOveral}}
%
	{\label{formulation:leader}}
%
	{}
%
%
	\addConstraint
	{\sum_{a \in \delta^+(i,bus)} f(a) z_a - \sum_{a \in \delta^-(i,bus)} f(a) z_a}
	{= 0 \label{eq:masterFlowConservation}}
	{\forall i \in V_H,}
	\addConstraint
	{\sum_{a \in A^{bus} \vert i(a) = i, j(a) = j} z_a}
	{\le 1 \label{eq:masterOneFrequency}}
	{\forall i, j \in V_H,}
	\addConstraint
	{x^t}
	{= \mathcal{C}^t(l^t) \quad \label{eq:userChoiceModel}}
	{\forall t \in T',}
	\addConstraint
	{x^t}
	{\in \mathbb{B} \label{eq:masterAdoption}}
	{\forall t \in T',}
	\addConstraint
	{z_a}
	{\in \mathbb{B} \label{eq:masterDesign}}
	{\forall a \in A^{bus},}
\end{mini!}%
\begin{center}
where $g^t$ and $l^t$ are the optimal values to the follower problem:
\end{center}%
\begin{customopti!}{\lexmin}
%
	{\boldsymbol{y}^t}
%
	{\left\langle \sum_{a \in A} \gamma_a^t y_{a}^t ~~,~ \sum_{a \in A} (\tau_a + \omega_a) y_a^t \right\rangle, \label{eq:objectivePassenger}}
%
	{\label{formulation:follower}}
%
	{\left\langle g^t, l^t\right\rangle = }
%
%
	\addConstraint
	{\sum_{{a} \in {\delta}^+({i})} y_{{a}}^t - \sum_{{a} \in {\delta}^-({i})} y_{{a}}^t}
	{= \begin{cases} 1 &\textrm{ if } {i} = o(t) \\ -1 &\textrm{ if } {i} = d(t) \\ 0 &\textrm{ else } \end{cases} \label{eq:extsubFlowConservation} \quad}
	{\forall {i} \in {V},}
	\addConstraint
	{y_{{a}}^t}
	{\le z_a \label{eq:extsubArcCapacity}}
	{\forall {a}  \in {A^{bus}},}
	\addConstraint
	{\sum_{a\in A} y_a^t}
	{\le K, \label{eq:extsubTransferLimit}}
	{}
	\addConstraint
	{y_{{a}}^t}
	{\in \mathbb{B} \label{eq:extsubVariables}}
	{\forall {a} \in {A}.}
\end{customopti!}
\caption{Bilevel Optimization Model for ODMTS Design with Adoptions}
\end{figure}

\paragraph{Solution Method}
Solving the bilevel model~\eqref{formulation:leader}-\eqref{formulation:follower} is computationally challenging for large-scale instances.
For this reason, this paper uses the iterative algorithm \texttt{arc-S2} by \citet{GuanEtAl2022-HeuristicAlgorithmsIntegrating} to approximate the optimal solution.
This algorithm splits the original bilevel framework into two components: (i) a regular ODMTS network design problem without considering latent demand and (ii) a choice model to analyze the adoption behavior of the potential riders.
The regular ODMTS network design problems are solved with the Benders-decomposition algorithm by \citet{auad2021resiliency}.
After the initial design of the ODMTS, a subset of the bus arcs is permanently fixed in the optimization, the behavior of the potential riders is analyzed, and the active set of riders is updated.
The ODMTS then undergoes a redesign based on the new set of riders, and this process is iterated until convergence.
The output is an ODMTS network: It defines bus and rail connections and their frequencies, and it defines which areas are served by on-demand shuttles.

\subsection{Including Dedicated Bus Lanes}
\label{subsec:dedicated_lane_ODMTS}

It only remains to include DBLs into the ODMTS network design.  Recall
that the travel time for congestion scenario $k$ is obtained by
scaling the travel time basis with the appropriate congestion factor.
For buses on DBLs, this travel time is simply replaced by the
free-flow travel time, after which the same optimization methods can
be used.  The travel time basis is not necessarily the free-flow
travel times: the baseline travel times are thus taken from
\citet{OpenStreetMap2021}.  It should be noted that the congestion
scenario and the availability of DBLs are both inputs to the network
design problem, such that the ODMTS delivers the best possible
performance under these circumstances.

The approach to include DBLs in the network is flexible in the sense that any road or set of roads can be replaced with DBLs for analysis.
The model is also not prescriptive in how the DBLs are realized, e.g., by using an existing lane or by adding a new lane.
This paper focuses on the impact on the public transit system, but an interesting direction for future research could be to study the trade-off with car congestion when a regular lane is converted to a DBL.
On the side of ODMTS, previous studies have shown that ODMTS generally reduce the total number of cars on the road \citep{auad2021resiliency,Vanhentenryck2023reach}.
Furthermore, increased traffic on the car lanes favors ODMTS, as the DBLs provide an even larger benefit.

\section{Case Study Gwinnett and Atlanta}
\label{sec:current}
\begin{figure}[p]
	\centering
	\includegraphics[width=0.7\linewidth]{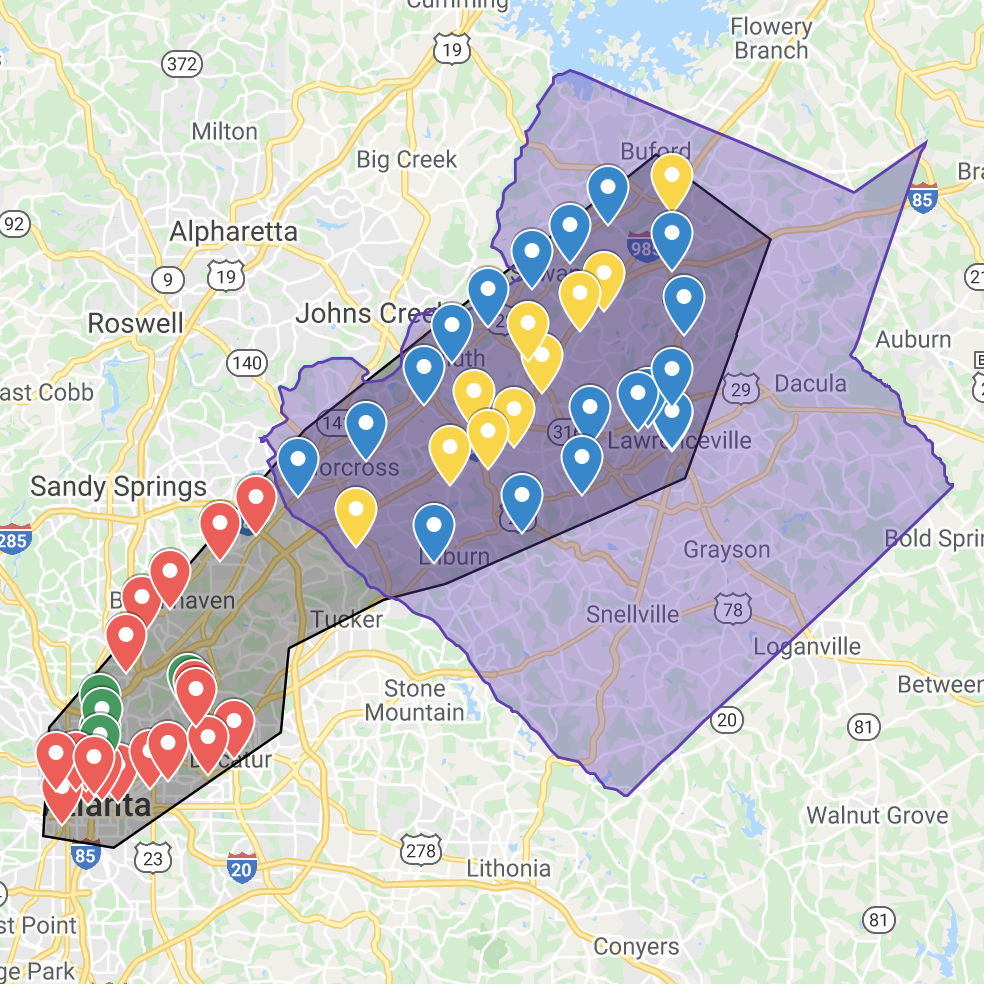}
	\caption{Case Study Area (Gray), Gwinnett County (Purple), and Potential Hub Locations (Markers)}
	\label{fig:study_area_hubs}
\end{figure}

The case study in this paper focuses on Interstate 85 (I-85) from Gwinnett County to the city of Atlanta, both in the Metro Atlanta Area.
The area is notorious for traffic, and Atlanta consistently ranks among the most congested cities in the US: After a temporary reduction in traffic during the COVID-19 pandemic, the city has jumped back to the 10th place on the list in 2021 \citep{Pishue2021-GlobalTrafficScorecard}.
Gwinnett County has almost a million residents, many of which work in Atlanta and drive the I-85 for their commute.
Gwinnett County recognizes the integral role of transit in the transportation system and published a transit development plan in 2018 \citep{GwinnettCounty2018-ComprehensiveTransitDevelopment}.
In 2019, Gwinnett County and the Georgia Department of Transportation (GDOT) presented a study of I-85 traffic and strategies to alleviate the extreme congestion \citep{Wickert2020-IsThereFix}.
However, the Atlanta Journal Constitution reports that Gwinnett County and GDOT officials are still seeking solutions to reduce congestion on I-85 \citep{Wickert2021-GwinnettGeorgiaDot}. This highly congested area serves as an excellent case to study the extent to which DBLs can improve passenger adoption in ODMTS.
Figure~\ref{fig:study_area_hubs} visualizes the case study area and its relation to Gwinnett County.
The area is specifically chosen to encompass I-85 and the immediate surrounding areas, and includes the current commuter routes operated by Gwinnett County Transit.
The markers represent strategic locations that will be used as potential hubs to design the ODMTS.
To clearly demonstrate the effect of congestion, the case study focuses on the morning rush hour from 7am to 9am and the associated southbound traffic from Gwinnett to Atlanta.

\paragraph{Current System}
\label{subsec:gwinnett-county-transit}

\begin{figure}[p]
	\centering
	\includegraphics[width=0.7\linewidth]{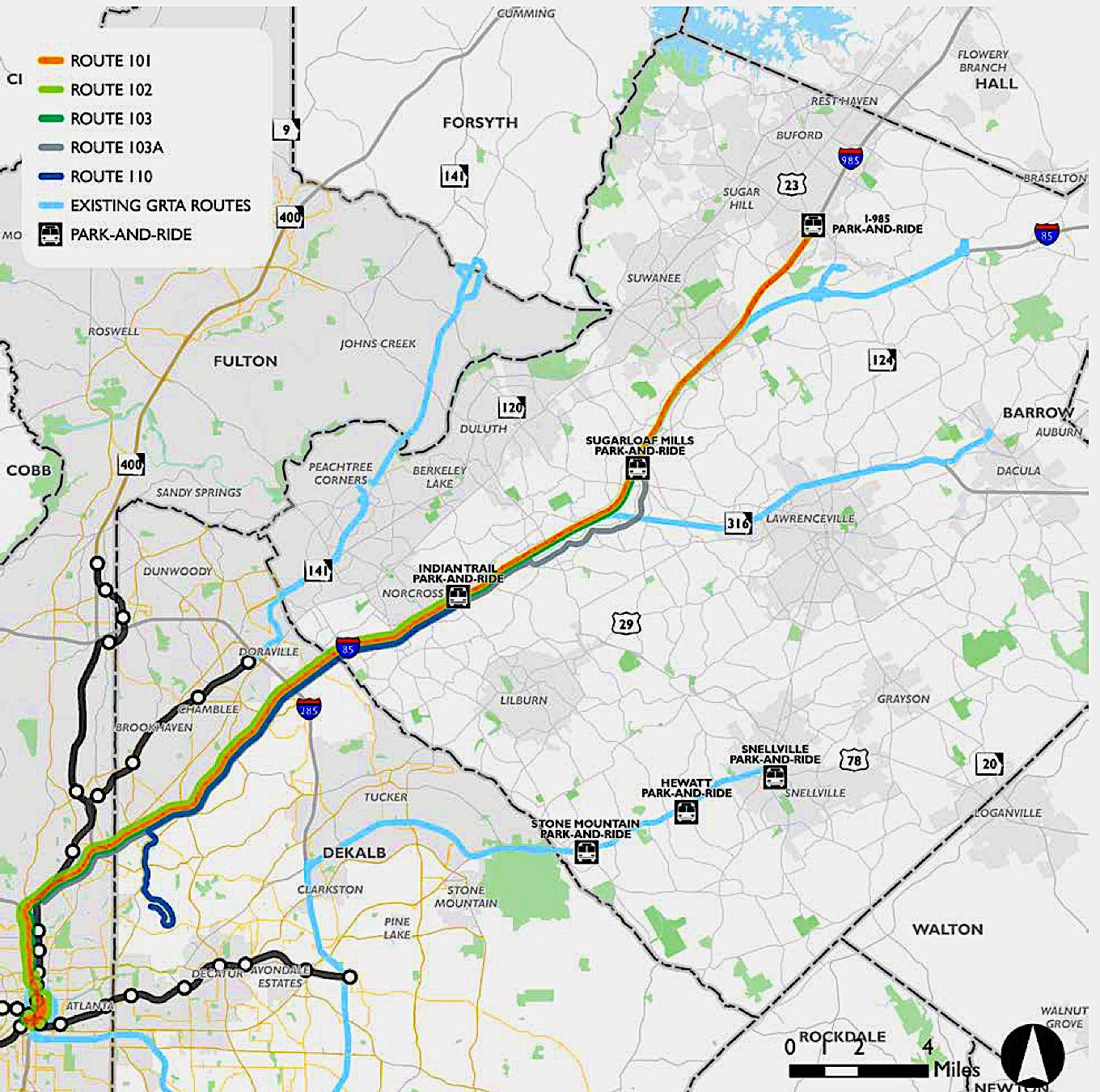}
	\caption{GCT Express Routes}
	\label{fig:express-routes}
\end{figure}

The bus transit system in Gwinnett County is operated by Gwinnett
County Transit (GCT).  GCT operates seven local bus routes within
Gwinnett County and five commuter routes that connect to DeKalb County
(Emory University) and Fulton County (Midtown Atlanta and Downtown
Atlanta).  The commuter routes are shown by
Figure~\ref{fig:express-routes}, and it can be seen that all routes
use the I-85.  In Atlanta, the buses stop at destinations in Emory
University, Midtown, and Downtown: they connect to the MARTA transit
system that is operated by the Metropolitan Atlanta Rapid Transit
Authority (MARTA), e.g., at the Civic Center station.  Where
available, transit vehicles make use of High Occupancy Toll (HOT)
lanes on the I-85 \citep{GDPS2022-I85Express}.  These lanes are shared
with three or more person carpools and with drivers who pay toll,
among other traffic.  To isolate the effect of DBLs, this study
compares the situation without HOT lanes to a new situation where the
extra lane is dedicated to buses only.

\paragraph{Existing Ridership}
This paper uses real transit data to generate a realistic set of trips
that represents the existing ridership. Historical ridership is
obtained from the Automated Fare Collection (AFC) system of GCT.
Riders use a transit card known as the \emph{Breeze Card} to tap onto
the buses.  For every tap, the AFC system records card id, time of
transaction, type of transaction, and which reader was tapped.  This
information is combined with the Automated Vehicle Location (AVL)
system that tracks the buses, and the Automated Passenger Counter
(APC) system that counts boardings and alightings at every stop.
Data was collected for the 7am-9am morning peak for April 16-19, 2018.
Trip chaining techniques from
\citet{BarryEtAl2002-OriginDestinationEstimation} are used to estimate
OD pairs.  Finally, the trips are sampled from the dataset to obtain a
set of 898 passengers that represent the existing ridership.  Together
these passengers generate 648 unique OD pairs.
The GCT transit data is the most representative and fine-grained data available for the existing ridership, but it does not track riders beyond their current boarding and alighting stops.
Existing riders may therefore see an additional improvement if they switch to ODMTS for their full door-to-door commute, eliminating the time it takes to get to and from the current transit stops.

\paragraph{Potential Riders}

In addition to existing riders, a set of potential riders is
generated: they may choose to adopt the system if the travel time is
sufficiently short.  The set of potential riders is based on simulated
travelers in the Metro Atlanta Area provided by the
\citet{ARC2022-ActivityBasedModeling}.  This data was generated by
an activity-based model, calibrated with data from the 2011 Regional
Houshold Travel Survey, and projected into 2020.
This dataset does not include the impact of the COVID-19 pandemic, but \citet{auad2021resiliency} have demonstrated that ODMTS are generally resilient to changes in demand.
To create a set of
potential transit users, travelers are selected from those who
\begin{itemize}
	\item commute alone by car from Gwinnett County;
	\item depart during the 7am-9am morning peak and stay within the case study area; and
	\item commute to the Emory University area, Midtown Atlanta, Downtown Atlanta, or stay within Gwinnett.
\end{itemize}
This results in 33,769 potential riders: 31,295 local riders (who
commute within Gwinnett) and 2,474 non-local riders who work outside
of the county.  The trip origins, which are provided at the level of
traffic analysis zones, are distributed over smaller census blocks
based on population counts.  Similarly, the trip destinations are
distributed over points of interest based on size, resulting in 27,202
unique OD pairs.
Note that potential ridership is based on simulated demand from the survey, while existing ridership is based on the current actual demand.

\paragraph{Experimental Settings}

The congestion scenarios are created according to
Section~\ref{subsec:congestion-modeling}, using Directions API data
for 8am on Wednesday March 22, 2022.  The MARTA rail system is fixed,
while the GCT buses are redesigned.  The ODMTS provides service to
1733 conveniently-placed \emph{virtual stops} that span the case study
area (chosen from existing transit stops, census block centers, and
points of interest), and the OD pairs are mapped to nearby virtual
stops.  Figure~\ref{fig:study_area_hubs} shows the potential hub
locations that are used for the case study.  These include hubs near
I-85 (yellow), hubs in Gwinnett (blue), MARTA rail stations (red), and
hubs at or near current GCT commuter bus stops in Atlanta (green).

Table~\ref{tab:odmts-paramters} shows the parameter values used for
the ODMTS design, with $\alpha$, $c^{bus}$ and $\bar{c}^{shuttle}$
taken from \citet{auad2021resiliency}.
The bus cost of \$72.15 per hour is based on data from \citet{FTA2018-NationalTransitDatabase} and \citet{Dickens2020-PublicTransportationVehicle} and consists of salaries and wages (\$39.55 per hour), maintenance (\$19.17 per hour), and vehicle depreciation (\$13.43 per hour).
The shuttle cost per mile is the estimated total cost to serve passengers with a real-time ride-sharing system.
The suggested shuttle cost per mile is multiplied by the congestion factor to capture the increase in
cost due to traffic.
Ticket prices are set to follow GCT.
Buses are added at a frequency of 10 buses per hour to maintain the notion of high-frequency connections. It is worth noting that because buses only serve the busy corridors, ODMTS can operate them at higher frequencies compared to typical routes.
\citet{auad2021resiliency} do not explicitly
model bus lines, and assume that every arc corresponds to a transfer
that induces wait time.  As bus lines from Gwinnett to Atlanta play an
important role, this paper removes the limit on how many arcs riders
can use ($K=\infty$) and adds waiting times in a post-processing step,
rather than at every stop.  Waiting for a bus or train is assumed to
be five minutes, which is a conservative estimate based on the three-minute expected waiting time and a two-minute buffer for transfers or to get into a station. Shuttles are assumed to be readily available, as discussed in Section~\ref{subsec:odmts-design}.
The set of adopting riders is updated accordingly in post-processing.

The ODMTS network design problem with adoptions is solved with the
iterative algorithm \texttt{arc-S2} by
\citet{GuanEtAl2022-HeuristicAlgorithmsIntegrating}, which is
implemented in Python 3.7.  In every iteration, the algorithm calls
the Benders-decomposition algorithm by \citet{auad2021resiliency} to
solve an ODMTS network design problem for fixed demand.  The
Benders-decomposition algorithm is implemented in C++ and uses CPLEX
12.9 to solve (mixed-integer) linear programs.

\begin{table}[!t]
	\centering
	\begin{tabular}{ll}
		\toprule
		Parameter & Value \\
		\midrule
		$\alpha$ & $0.1078$ (value time at \$7.25/hour to match U.S. federal minimum wage). \\
		$c^{bus}$ & \$72.15 per hour. \\
		$\bar{c}^{shuttle}$ & Basis of \$1/mile multiplied by congestion scaling factor. \\
		$c^{ticket}$ & \$3.50/trip for local and \$5/trip for non-local.\\
		$F^{bus}$ & \{10\} buses/hour.\\
		$K$ & $\infty$ (no transfer limit).\\
		$\rho$ & 1.5 adoption factor.\\
		\bottomrule
	\end{tabular}
	\caption{ODMTS Parameters for Case Study}
	\label{tab:odmts-paramters}
\end{table}

\section{Impact of Congestion on ODMTS}
\label{sec:baseline}

This section studies the impact of congestion on ODMTS.  For each of
the congestion scenarios, an ODMTS is designed to deliver optimal
performance in that scenario, and it is analyzed how these systems
differ in service quality and cost.  This analysis also serves as
the baseline for studying the impact of DBLs in
Section~\ref{sec:impact-of-DBLs}.

\paragraph{Network Designs}

\begin{figure}[!t]
	\centering
	\begin{subfigure}{.47\columnwidth}
		\centering
		\includegraphics[width=\textwidth]{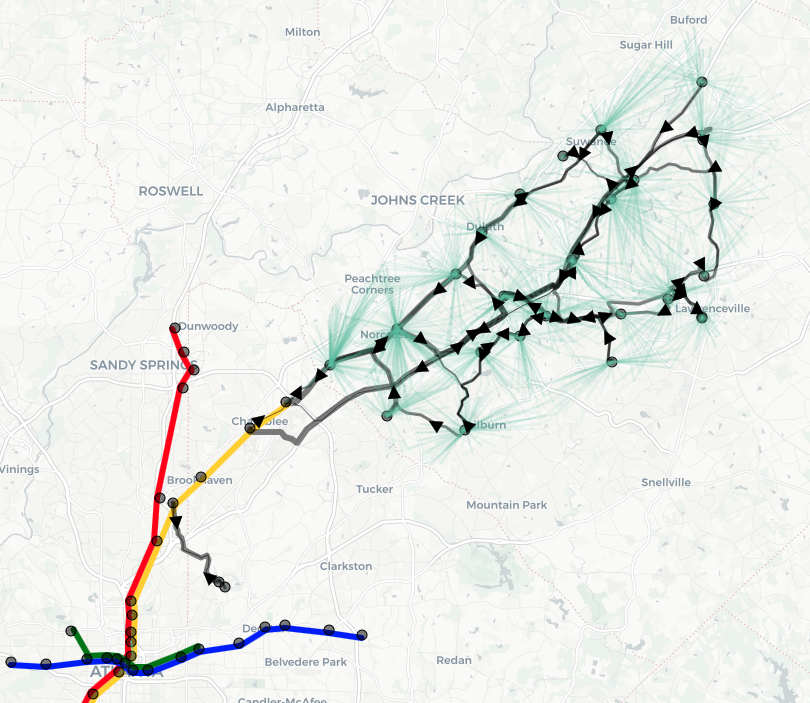}
		\caption{Expected Scenario}
		\label{fig:adoption-map-expected-ndl}
	\end{subfigure}%
	\hspace{0.05\textwidth}
	\begin{subfigure}{.47\columnwidth}
		\centering
		\includegraphics[width=\textwidth]{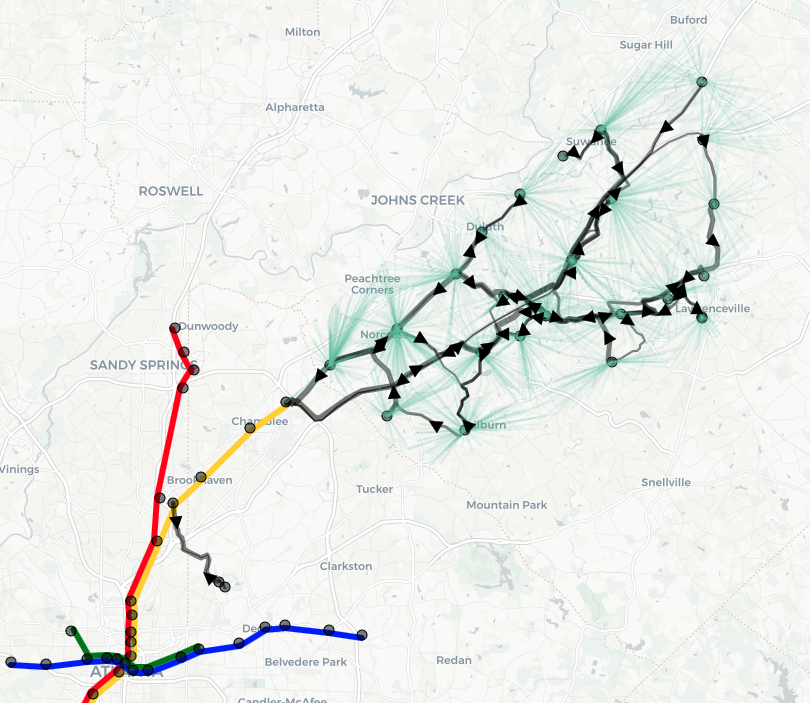}
		\caption{50-50 Scenario}
		\label{fig:adoption-map-50-50-ndl}
	\end{subfigure}\\
	\vspace{\baselineskip}
	\begin{subfigure}{.47\columnwidth}
		\centering
		\includegraphics[width=\textwidth]{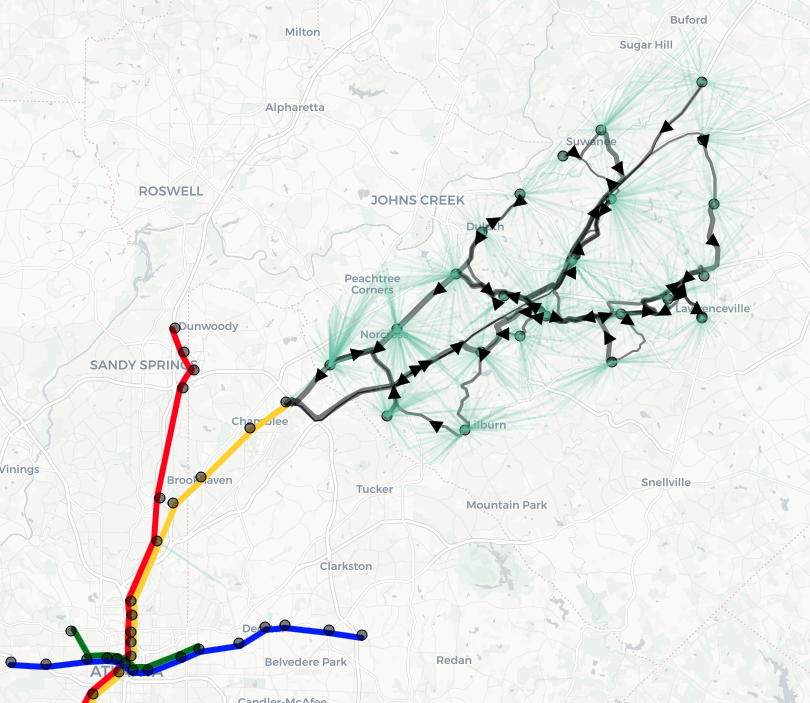}
		\caption{Pessimistic Scenario}
		\label{fig:adoption-map-pessimistic-ndl}
	\end{subfigure}\\
	\caption{ODMTS Designs without DBLs under Three Congestion Scenarios}
	\label{fig:adoption-maps-ndl}
\end{figure}

Figure~\ref{fig:adoption-maps-ndl} presents the ODMTS network designs,
optimized under the different congestion scenarios.  The lines in the
Southwest represent the Red, Gold, Blue, and Green lines of the MARTA
rail system.  The lines with arrows (gray) show the bus arcs that are
opened by the optimization model, where the thickness of the line
corresponds to the number of travelers that use this arc.  Finally,
the thin lines (green) indicate which connections are served by
on-demand shuttles.  It is striking that none of the designs feature a
bus connection on the I-85 from Gwinnett to Atlanta.  This corridor is
so congested that the ODMTS transfers passengers to rail at Chamblee
and Doraville stations rather than provide a direct bus connection to
the city.  In the expected and 50-50 scenarios, there is a bus service
that connects the rail to Emory University. Under the pessimistic
scenario, this bus also disappears in favor of serving a small number
of people by shuttles.  The overall designs look similar for the
different scenarios, but the ODMTS increases its reliance on shuttles
as congestion increases: the total shuttle distance increases from 62k
miles for the expected scenario, to 66k miles for the 50-50 scenario,
and 72k miles for the pessimistic scenario.

\paragraph{Travel Time and Adoption}

\begin{figure}[t]
    \centering
    \begin{minipage}{0.47\textwidth}
    	\centering
	    \includegraphics[width=\linewidth,trim={1cm 0 0 0},clip]{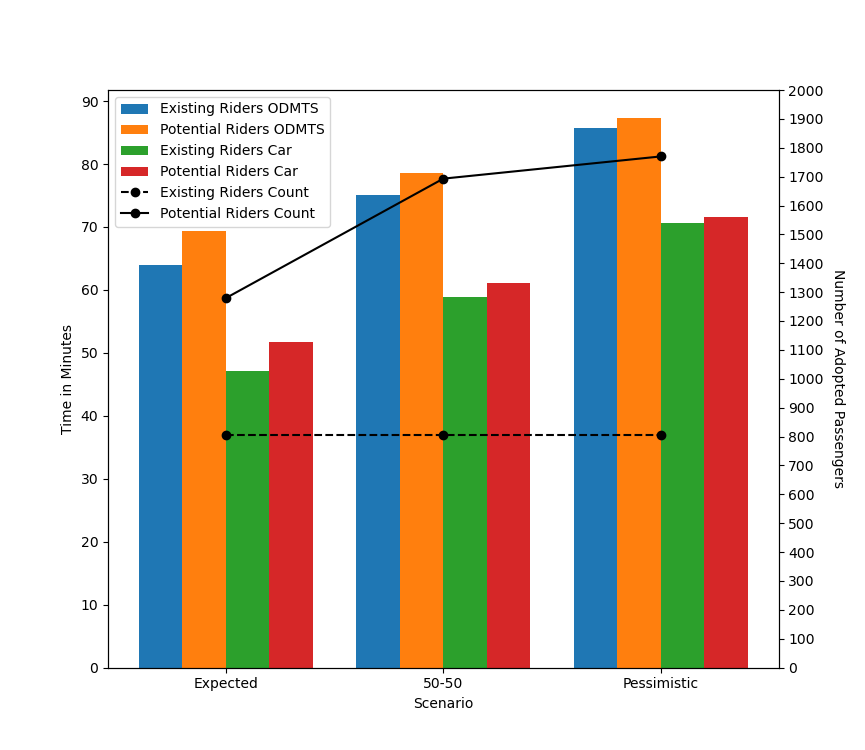}
        \captionof{figure}{Baseline Travel Times (left axis) and Counts (right axis) for Non-local Riders}
	    \label{fig:odmts-in-gwinnett-travel-times}
	\end{minipage}%
	\hspace{0.05\textwidth}
    \begin{minipage}{0.47\textwidth}
		\centering
		\includegraphics[width=\linewidth]{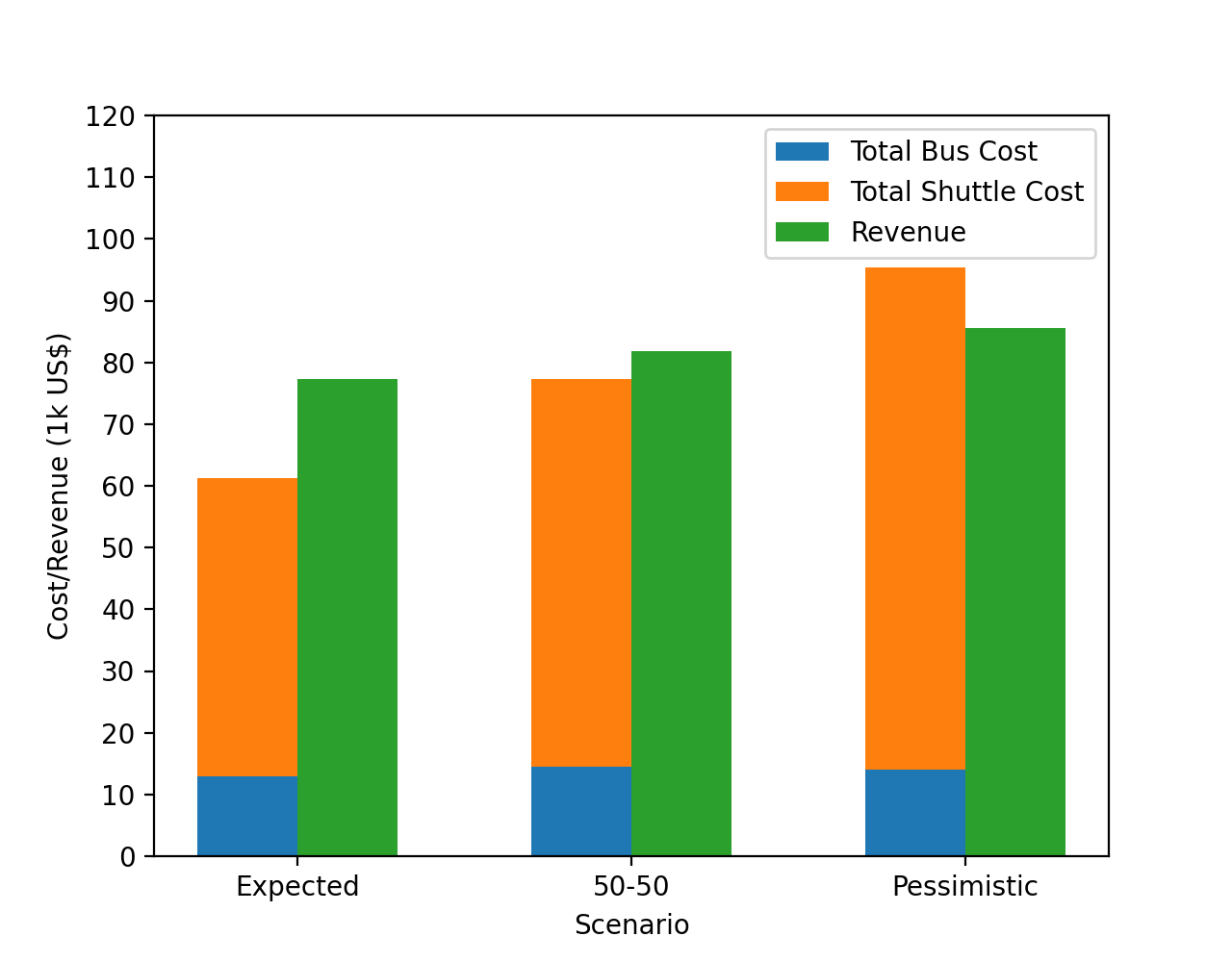}
		\caption{Baseline ODMTS System Costs}
		\label{fig:odmts-in-gwinnett-transit-costs}
\end{minipage}
\end{figure}

Figure~\ref{fig:odmts-in-gwinnett-travel-times} zooms in on the travel
time and adoption of non-local travelers, who are most likely to be
affected by congestion.  Congestion has a major effect on the existing
riders: travel time increases from 64 minutes under the expected
scenario to 86 minutes under the pessimistic scenario.  While the
number of existing riders is fixed, potential riders adopt the ODMTS
if the travel time is not too long compared to driving.  Surprisingly,
Figure~\ref{fig:odmts-in-gwinnett-travel-times} shows an increase in
adoption when congestion increases, despite the increase in travel
time.  This is due to the fact that travel time by car is also
negatively impacted, making transit a more favorable option.  There is
52\% adoption of potential non-local riders under the expected
scenario (1279 people), which increases to 72\% under the pessimistic
scenario (1770 people).  Local riders are less affected by congestion,
but show similar trends.  Travel time for existing riders increases
from 23 minutes under the expected scenario to 33 minutes under the
pessimistic scenario, while the adoption rate increases from 61\% to
66\% of the potential riders (from 19,013 to 20,646 riders).  The high
adoption rates suggest that ODMTS may significantly improve access to
transit and may compete with commuting by car for a variety of
congestion levels.

\begin{figure}[t]
	\centering
	\begin{subfigure}{0.47\textwidth}
		\centering
		\includegraphics[width=\linewidth]{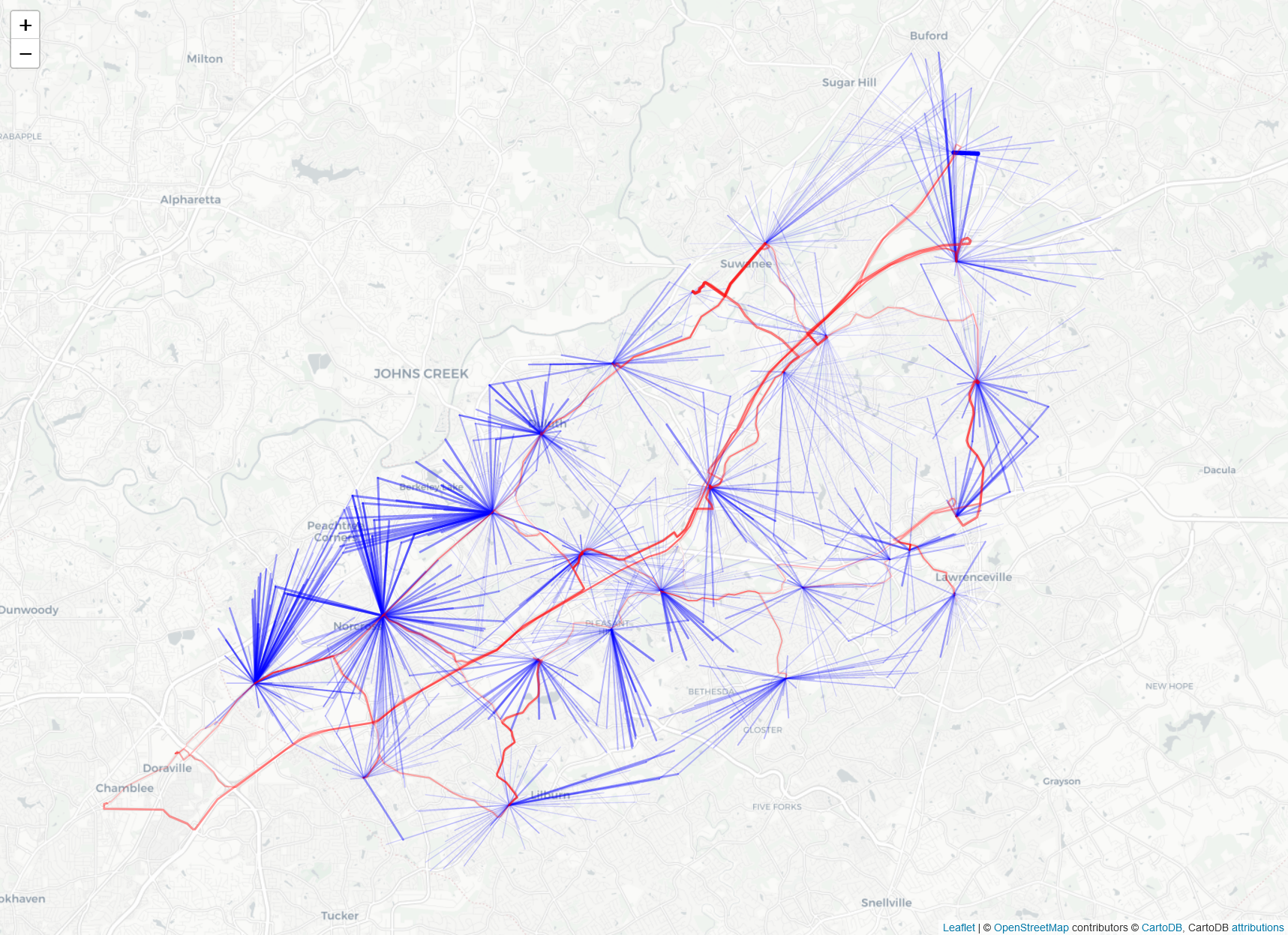}
		\caption{Expected Scenario}
		\label{fig:odmts-in-gwinnett-scaling-factors:expected}
	\end{subfigure}%
	\hspace{0.05\textwidth}
	\begin{subfigure}{0.47\textwidth}
		\centering
		\includegraphics[width=\linewidth]{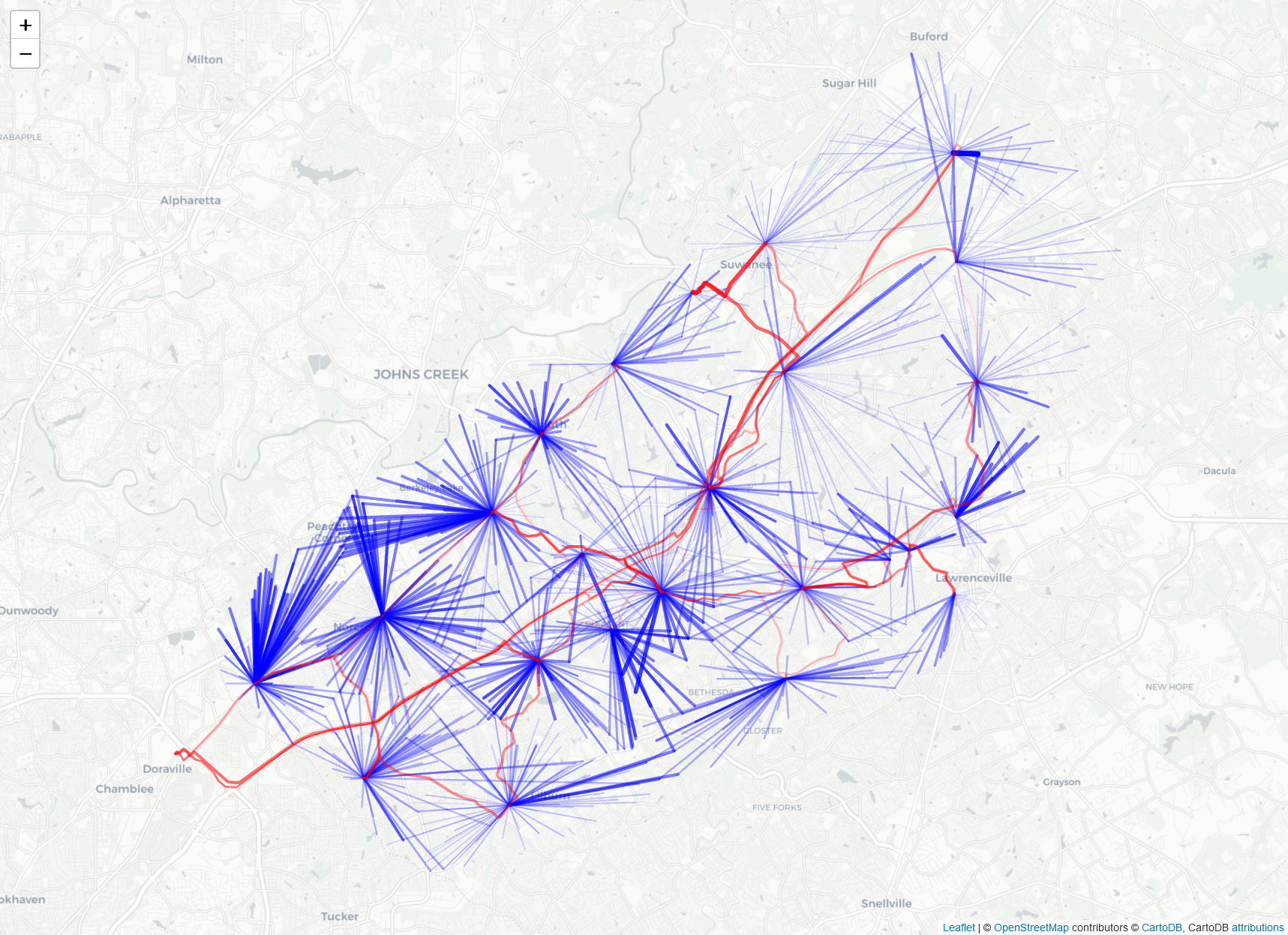}
		\caption{Pessimistic Scenario}
		\label{fig:odmts-in-gwinnett-scaling-factors:pessimistic}
	\end{subfigure}%
	\caption{Baseline ODMTS Congestion Scaling Factors of Used Roads (darker is higher)}
	\label{fig:odmts-in-gwinnett-scaling-factors}
\end{figure}

\paragraph{Congestion}
To get a better understanding of how congestion affects system
operations, Figure~\ref{fig:odmts-in-gwinnett-scaling-factors}
displays the scaling factors for the roads used by the ODMTS.  For the
expected scenario
(Figure~\ref{fig:odmts-in-gwinnett-scaling-factors:expected}) the
average scaling factors for roads used by shuttles and buses are
$R^{shuttle}=0.78$ and $R^{bus}=0.89$, respectively.
These are averages over the arcs, weighted by distance.
In the pessimistic scenario
(Figure~\ref{fig:odmts-in-gwinnett-scaling-factors:pessimistic}) these
scaling factors go up significantly to $R^{shuttle}=1.13$ and
$R^{bus}=1.17$.
It is important to note that these networks are
designed while taking congestion into account.  That is, the optimal
designs avoid the extremely congested I-85 and distribute the
remaining burden of traffic about equally over the shuttles and the
buses.  As congestion increases, substituting buses for shuttles makes
this possible.

\paragraph{System Cost} Figure~\ref{fig:odmts-in-gwinnett-transit-costs} details the cost and revenue of the ODMTS under the different scenarios.
Ridership increases with congestion, which leads to a small increase in revenue.
The costs increase at the same time, mainly because of the increased travel time for buses and shuttles.
Shuttle costs play the most prominent role, making up 79\% to 85\% of the total cost.
The bus costs stay relative constant because the ODMTS increases its reliance on shuttles as congestion increases.
Figure~\ref{fig:odmts-in-gwinnett-transit-costs} also indicates that the ODMTS is profitable for the first two scenarios, but has a net cost of \$0.42 per person in the pessimistic scenario.
However, note that traditional transit systems are often heavily subsidized, where the ODMTS almost breaks even while serving more than 20k additional passengers and providing better access to transit.
Furthermore, \citet{auad2021resiliency} show that a stricter budget can be enforced by reducing the number of on-demand shuttles at the cost of a relatively small increase in waiting time.

\section{Impact of Dedicated Bus Lanes}
\label{sec:impact-of-DBLs}

The previous section has shown that without DBLs, congestion has a
major impact on ODMTS travel times. This section studies to what
extent DBLs can mitigate these negative effects, even when shuttles
are still affected by traffic.  The impact of DBLs is assessed by
opening a DBL on I-85, reoptimizing the networks, and evaluating their
performance.  The ODMTS with DBLs is compared to the baseline in terms
of design, travel times, passenger adoption, and cost of the system.
Finally,
this section explores how adoption changes with the adoption
factor $\rho$.

\paragraph{Network Redesigns with DBLs}
\label{sssec:odmts-network}

\begin{figure}[!t]
	\centering
	\begin{subfigure}{.47\columnwidth}
		\centering
		\includegraphics[width=\linewidth]{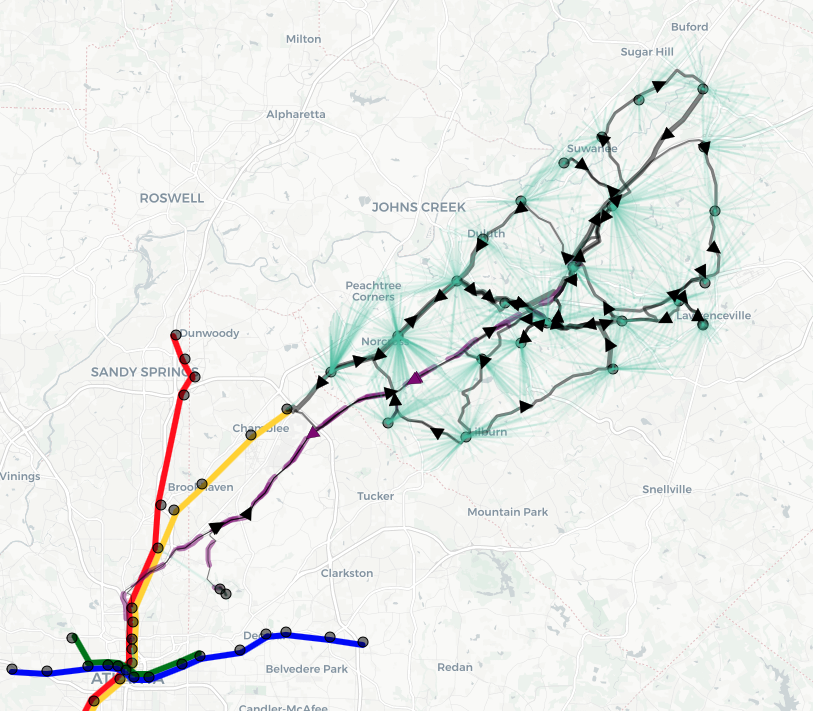}
		\caption{Expected Scenario}
		\label{fig:adoption-map-expected-dl}
	\end{subfigure}%
\hspace{0.05\textwidth}
	\begin{subfigure}{.47\columnwidth}
		\centering
		\includegraphics[width=\linewidth]{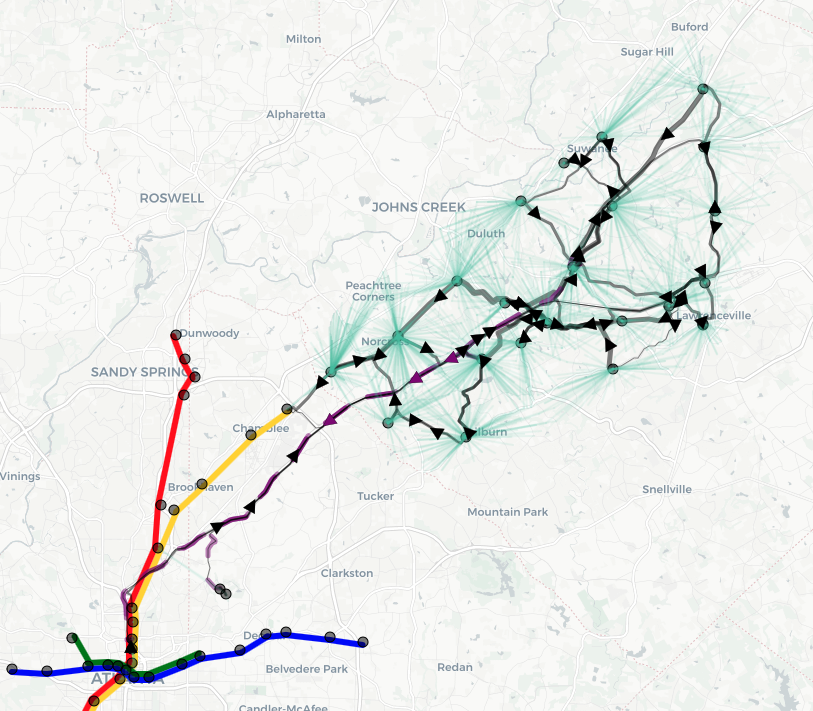}
		\caption{50-50 Scenario}
		\label{fig:adoption-map-50-50-dl}
	\end{subfigure}\\
	\vspace{\baselineskip}
	\begin{subfigure}{.47\columnwidth}
		\centering
		\includegraphics[width=\linewidth]{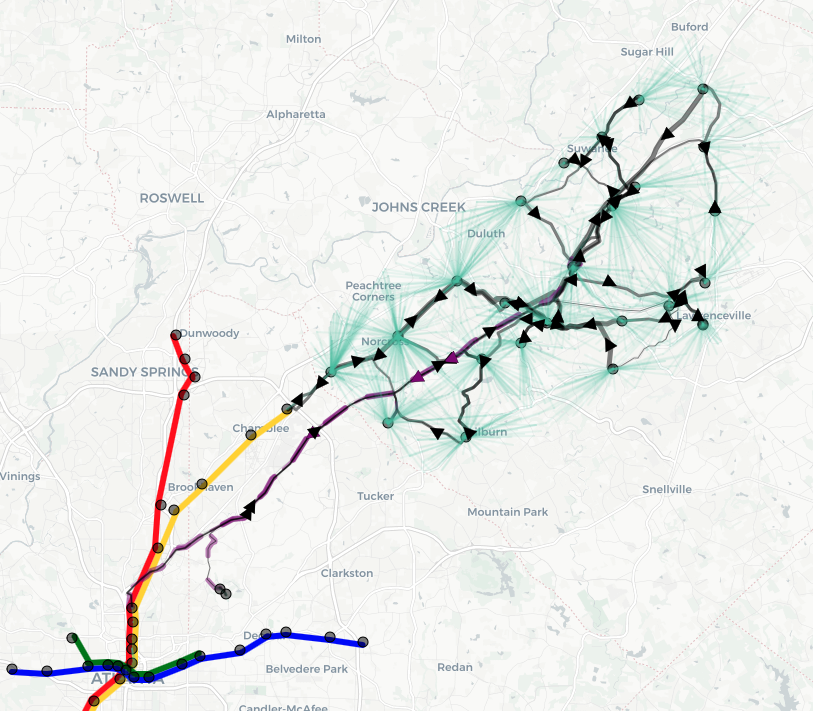}
		\caption{Pessimistic Scenario}
		\label{fig:adoption-map-pessimistic-dl}
	\end{subfigure}\\
	\caption{ODMTS Designs with DBLs under Three Congestion Scenarios}
	\label{fig:adoption-maps-dl}
\end{figure}

Figure~\ref{fig:adoption-maps-dl} presents the ODMTS network redesigns
that are optimized to leverage the DBLs (indicated by dashed lines).  As a major change from the
baseline ODMTS, the ODMTS with DBLs now include buses traveling south
on I-85 from Gwinnett to Atlanta for all scenarios.  This corridor was
too congested to be used in the baseline, but the thick dashed lines in
Figure~\ref{fig:adoption-maps-dl} indicate that riders now benefit
extensively from these new bus connections.  The local network within
Gwinnett is similar across the scenarios, and also similar to the
baseline. Recall that the baseline sees a substitution effect of buses
for shuttles when congestion increases. This effect largely disappears
as the important bus lanes are no longer affected by congestion
(shuttle mileages are 65k, 66k, and 67k for the three scenarios,
respectively).  Overall the new networks suggest a major improvement
for non-local travelers, and similar service for local travelers.

\paragraph{Travel Times and Adoption}

\begin{figure}[!t]
	\centering
	\resizebox{\textwidth}{!}{
		\begin{tabular}{llccccccccc}
			\toprule
			&       & \multicolumn{5}{c}{Ridership}         & \multicolumn{4}{c}{Average Travel Time (minutes)} \\
			\cmidrule(lr){3-7} \cmidrule(lr){8-11}
			&       & Existing & \multicolumn{4}{c}{Adoption} & Car   & \multicolumn{3}{c}{ODMTS} \\
			\cmidrule(lr){3-3} \cmidrule(lr){4-7} \cmidrule(lr){8-8} \cmidrule(lr){9-11}
			Scenario & DBL   & Count & Count & Rate & PreDBL & PostDBL & Existing & Existing & PreDBL & PostDBL \\
			\midrule
			\multirow{2}[0]{*}{Expected} & No    & 806   & 1,279 & 52\%  & -     & -     & 47    & 64    & 70    & 70 \\
			& Yes   & 806   & 2,078 & 84\%  & 1,263 & 815   & 47    & 43    & 54    & 47 \\
			\midrule
			\multirow{2}[0]{*}{50-50} & No    & 806   & 1,693 & 68\%  & -     & -     & 59    & 75    & 79    & 85 \\
			& Yes   & 806   & 2,153 & 87\%  & 1,665 & 488   & 59    & 44    & 56    & 53 \\
			\midrule
			\multirow{2}[0]{*}{Pessimistic} & No    & 806   & 1,770 & 72\%  & -     & -     & 71    & 86    & 88    & 107 \\
			& Yes   & 806   & 2,307 & 93\%  & 1,752 & 555   & 71    & 52    & 61    & 60 \\
			\bottomrule
		\end{tabular}%
	}
	\captionof{table}{Adoption and Travel Times for Non-local Riders}
	\label{tab:adoption-traveltimes-nonlocal}%
	\hspace{\baselineskip}
	\resizebox{\textwidth}{!}{
		\begin{tabular}{llccccccccc}
			\toprule
			&       & \multicolumn{5}{c}{Ridership}         & \multicolumn{4}{c}{Average Travel Time (minutes)} \\
			\cmidrule(lr){3-7} \cmidrule(lr){8-11}
			&       & Existing & \multicolumn{4}{c}{Adoption} & Car   & \multicolumn{3}{c}{ODMTS} \\
			\cmidrule(lr){3-3} \cmidrule(lr){4-7} \cmidrule(lr){8-8} \cmidrule(lr){9-11}
			Scenario & DBL   & Count & Count & Rate & PreDBL & PostDBL & Existing & Existing & PreDBL & PostDBL \\
			\midrule
			\multirow{2}[0]{*}{Expected} & No    & 92    & 19,013 & 61\%  & -     & -     & 10    & 20    & 11    & 31 \\
			& Yes   & 92    & 19,091 & 61\%  & 17,942 & 1,149 & 10    & 18    & 11    & 19 \\
			\midrule
			\multirow{2}[0]{*}{50-50} & No    & 92    & 19,706 & 63\%  & -     & -     & 13    & 23    & 14    & 35 \\
			& Yes   & 92    & 19,266 & 62\%  & 18,121 & 1,145 & 13    & 22    & 14    & 21 \\
			\midrule
			\multirow{2}[0]{*}{Pessimistic} & No    & 92    & 20,646 & 66\%  & -     & -     & 15    & 28    & 16    & 48 \\
			& Yes   & 92    & 19,618 & 63\%  & 18,295 & 1,323 & 15    & 25    & 16    & 31 \\
			\bottomrule
		\end{tabular}%
	}
	\captionof{table}{Adoption and Travel Times for Local Riders}
	\label{tab:adoption-traveltimes-local}
\end{figure}

Tables~\ref{tab:adoption-traveltimes-nonlocal} and
\ref{tab:adoption-traveltimes-local} support the view that DBLs bring
significant benefits to non-local travelers while local service levels
remain high. The tables separate adopted riders into those who would
adopt even without DBLs (PreDBL) and those who would only switch to
ODMTS if DBLs are implemented (PostDBL).
Table~\ref{tab:adoption-traveltimes-nonlocal} shows huge improvements
for non-local travelers in terms of both travel time and adoption.
Without DBLs, existing riders are expected to spend 64 minutes in the
ODMTS, compared to 47 minutes by car.  With DBLs, the average ODMTS
travel time goes down to only 43 minutes and becomes faster than
driving.  This is possible because the buses are able to avoid
traffic.  PreDBL and PostDBL adopters also benefit from significantly
improved travel times compared to ODMTS without DBLs.
For example, PostDBL riders would have an average travel time of 107 minutes if they use ODMTS without DBLs, which reduces to only 60 minutes when DBLs are implemented.
In fact, \emph{DBLs are so effective
  that the travel times under the pessimistic scenario are better than
  the expected scenario without DBLs}.  This is strong evidence that
DBLs can mitigate the effects of congestion, even when shuttles are
still affected by traffic.

The lower travel times for non-local riders result in a significantly
increased adoption rate.  For the expected scenario, the increase is 32
percentage point (from 52\% to 84\%), and for the pessimistic scenario
it increases from 72\% to 93\%.  This suggests that, in a congested
environment, DBLs are able to attract almost all potential riders that
can benefit from the traffic-free lane.
Table~\ref{tab:mode-distributions} provides additional information on
the mode distribution of the non-local adopters.  Without DBLs, the
only viable option is to take a shuttle to a bus station, take the bus
to Chamblee or Doraville station and transfer to the rail.  With DBLs,
it becomes viable to directly take the bus to Atlanta rather than take
a detour by rail, and this attracts a significant number of additional
riders.  Using the rail is still popular but, rather than transferring
at Chamblee or Doraville station, riders get on the rail in Midtown
Atlanta to complete the last mile.

\begin{table}[!t]
	\centering
	\begin{tabular}{lcccccc}
		\toprule
		& \multicolumn{2}{c}{Expected}     & \multicolumn{2}{c}{50-50}         & \multicolumn{2}{c}{Pessimistic}   \\
		\cmidrule(lr){2-3} \cmidrule(lr){4-5} \cmidrule(lr){6-7}
		& No DBL & DBL & No DBL  & DBL & No DBL  & DBL \\
		\midrule
		Bus               & 0   & 1      & 0    & 1      & 0    & 1      \\
		Bus and Rail            & 4   & 5      & 6    & 4      & 5    & 4      \\
		Shuttle and Bus         & 0 & 674     & 0  & 754     & 0  & 876     \\
		Shuttle, Bus, and Rail & 1275 & 1398     & 1687  & 1394     & 1765 & 1426    \\
		\midrule
		Total       & 1279 & 2078    & 1693 & 2153    & 1770 & 2307\\
		\bottomrule
	\end{tabular}
	\caption{Mode Distribution for Non-local Adopted Potential Riders}
	\label{tab:mode-distributions}
\end{table}

Non-local riders receive the most benefit from DBLs, but Table~\ref{tab:adoption-traveltimes-local} shows that these improvements do not come at the expense of local transit users.
Existing local riders see a small improvement in average travel time (up to 3 minutes in the pessimistic scenario), and the travel times for PreDBL adopters are more or less constant.
Adoption of the ODMTS stays high, with over 60\% of potential riders adopting in all cases.
However, there is some decrease in adoption when DBLs are implemented in the 50-50 scenario (63\% to 62\%) and in the pessimistic scenario (66\% to 63\%).
On the other hand, the DBL redesign also brings in new riders by significantly lowering their average travel time by up to 17 minutes in the pessimistic case.
Overall the service level for the local travelers is comparable between the baseline and the DBL redesign, while the service level for non-local travelers improves significantly.

\paragraph{System Cost}
\label{sssec:cost-design}

\begin{figure}[th]
    \centering
    \includegraphics[width=.7\linewidth]{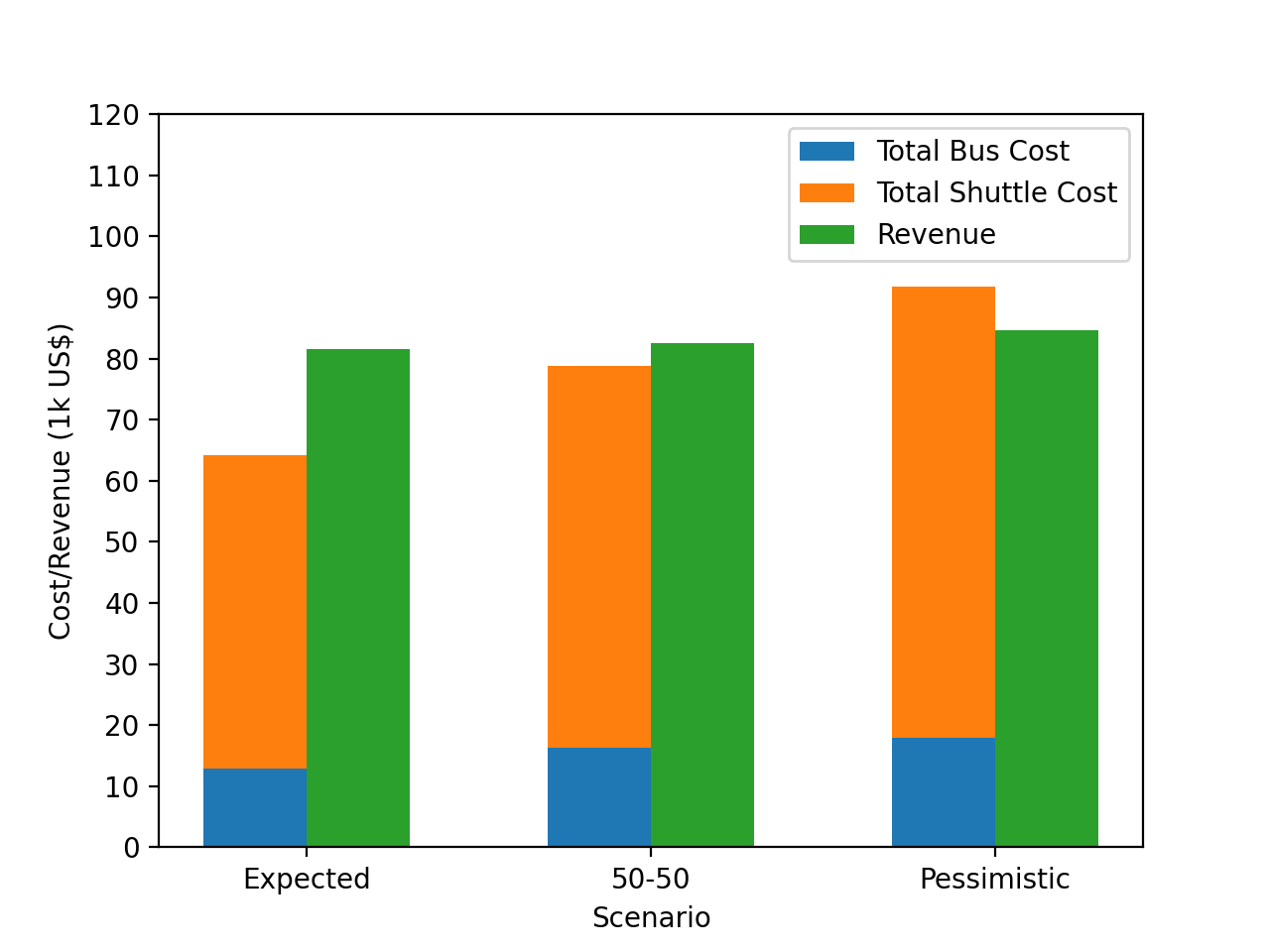}
    \caption{ODMTS with DBLs System Costs}
    \label{fig:odmts-redesign-transit-costs}
\end{figure}

Table~\ref{tab:cost-stats} shows how DBLs impact the cost of the
system.  As congestion increases, DBLs become more useful compared to
the baseline, and the table shows bigger investments in buses.  For
example, in the expected scenario, the redesign increases the shuttle
cost by \$3k while the bus cost remains the same.  In the pessimistic
scenario, on the other hand, investing \$4k into buses allows the
shuttle cost to be \emph{reduced} by \$7k, leading to a lower total
cost.  For the expected scenario and the 50-50 scenario, the increase
in revenue balances out the increase in cost, resulting in only \$0.03
difference in net profit per rider.  For the pessimistic scenario, the
lower system cost results in a savings of \$0.10 per passenger. As a
result, DBLs are effective in reducing travel times and increasing
adoption without negatively affecting system cost.  At the same time,
DBLs do not resolve the issue that system cost increases with
congestion, as depicted in Figure~\ref{fig:odmts-redesign-transit-costs}.  However, it is worth repeating that traditional transit
systems tend to be heavily subsidized, where the ODMTS almost breaks
even when serving 20k additional passengers.

\begin{table}[!t]
	\centering
	\resizebox{\textwidth}{!}{
		\begin{tabular}{llcccccc}
			\toprule
			&       &       & \multicolumn{3}{c}{Cost} &       &  \\
			\cmidrule(lr){4-6}
			Scenario & DBL   & Ridership & Total & Bus   & Shuttle & Revenue & Net Profit per Rider \\
			\midrule
			\multirow{2}[0]{*}{Expected} & No    & 20,292 & \$ 61k & \$ 13k & \$ 48k & \$ 77k & \$ 0.76 \\
			& Yes   & 21,169 & \$ 64k & \$ 13k & \$ 51k & \$ 82k & \$ 0.79 \\
			\midrule
			\multirow{2}[0]{*}{50-50} & No    & 21,399 & \$ 77k & \$ 14k & \$ 63k & \$ 82k & \$ 0.20 \\
			& Yes   & 21,419 & \$ 79k & \$ 16k & \$ 63k & \$ 83k & \$ 0.17 \\
			\midrule
			\multirow{2}[0]{*}{Pessimistic} & No    & 22,416 & \$ 95k & \$ 14k & \$ 81k & \$ 85k & -\$ 0.42 \\
			& Yes   & 21,925 & \$ 92k & \$ 18k & \$ 74k & \$ 85k & -\$ 0.32 \\
			\bottomrule
		\end{tabular}%
	}
	\caption{ODMTS System Costs}
	\label{tab:cost-stats}
\end{table}

\subsection{Sensitivity Analysis on Adoption Factor}
\label{sec:sensitivity_adoption}
\begin{figure}[!t]
	\centering
	\begin{subfigure}[h]{0.33\textwidth}
		\centering
		\includegraphics[width=\textwidth]{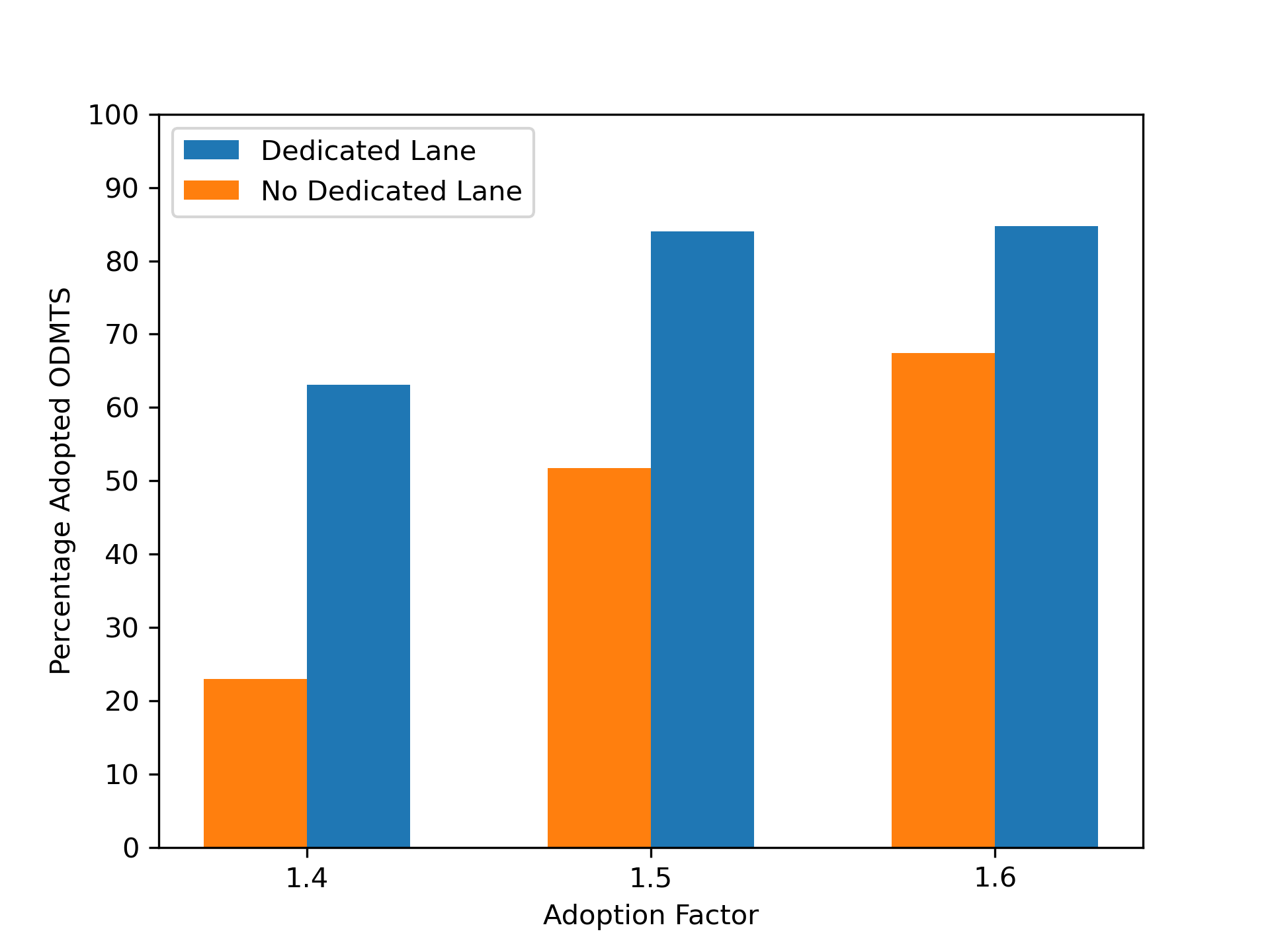}
		\caption{Expected Scenario}
		\label{fig:adoption-percent-expected}
	\end{subfigure}
	\begin{subfigure}[h]{0.33\textwidth}
		\centering
		\includegraphics[width=\textwidth]{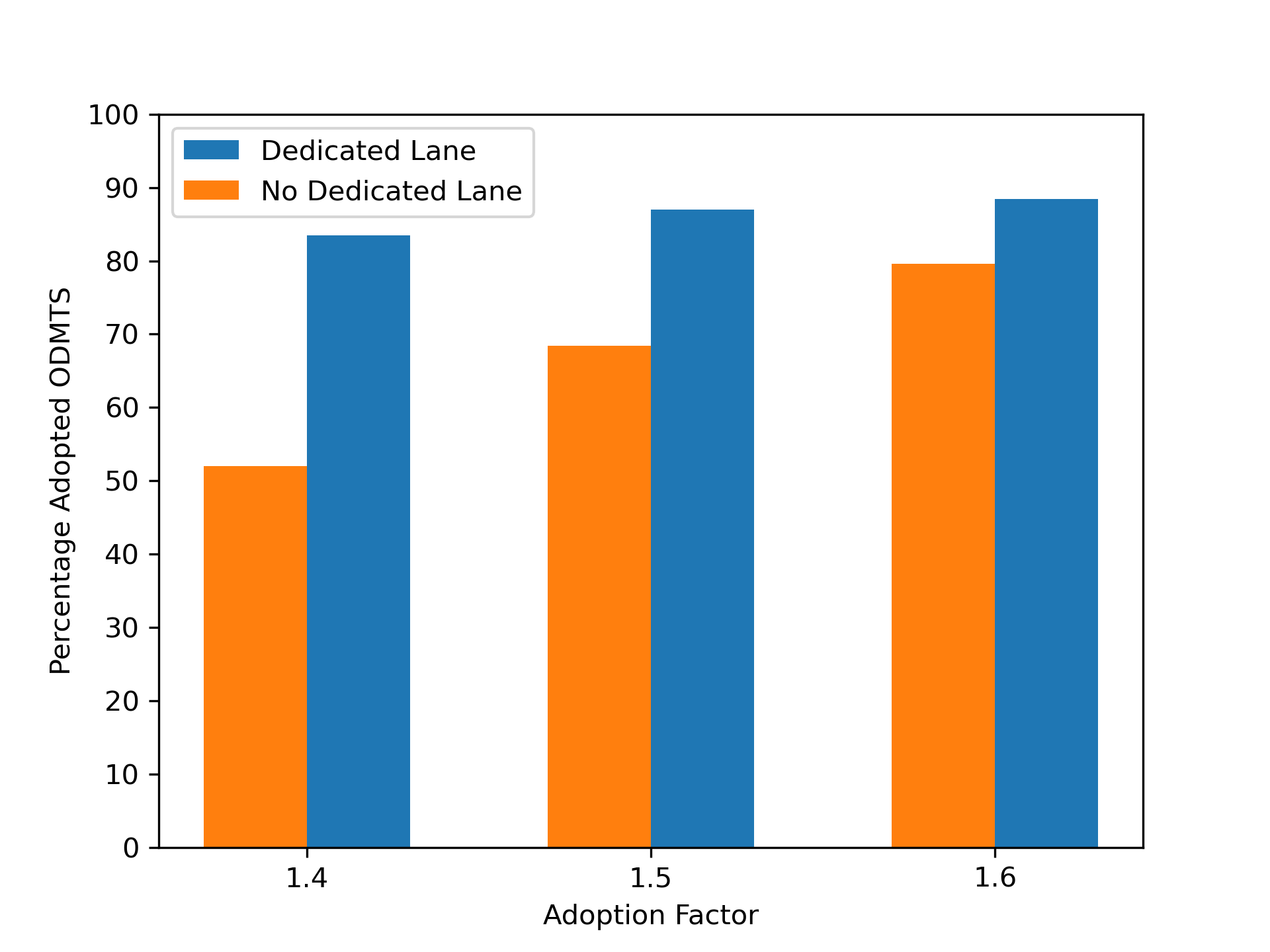}
		\caption{50-50 Scenario}
		\label{fig:adoption-percent-50-50}
	\end{subfigure}
	\begin{subfigure}[h]{0.33\textwidth}
		\centering
		\includegraphics[width=\textwidth]{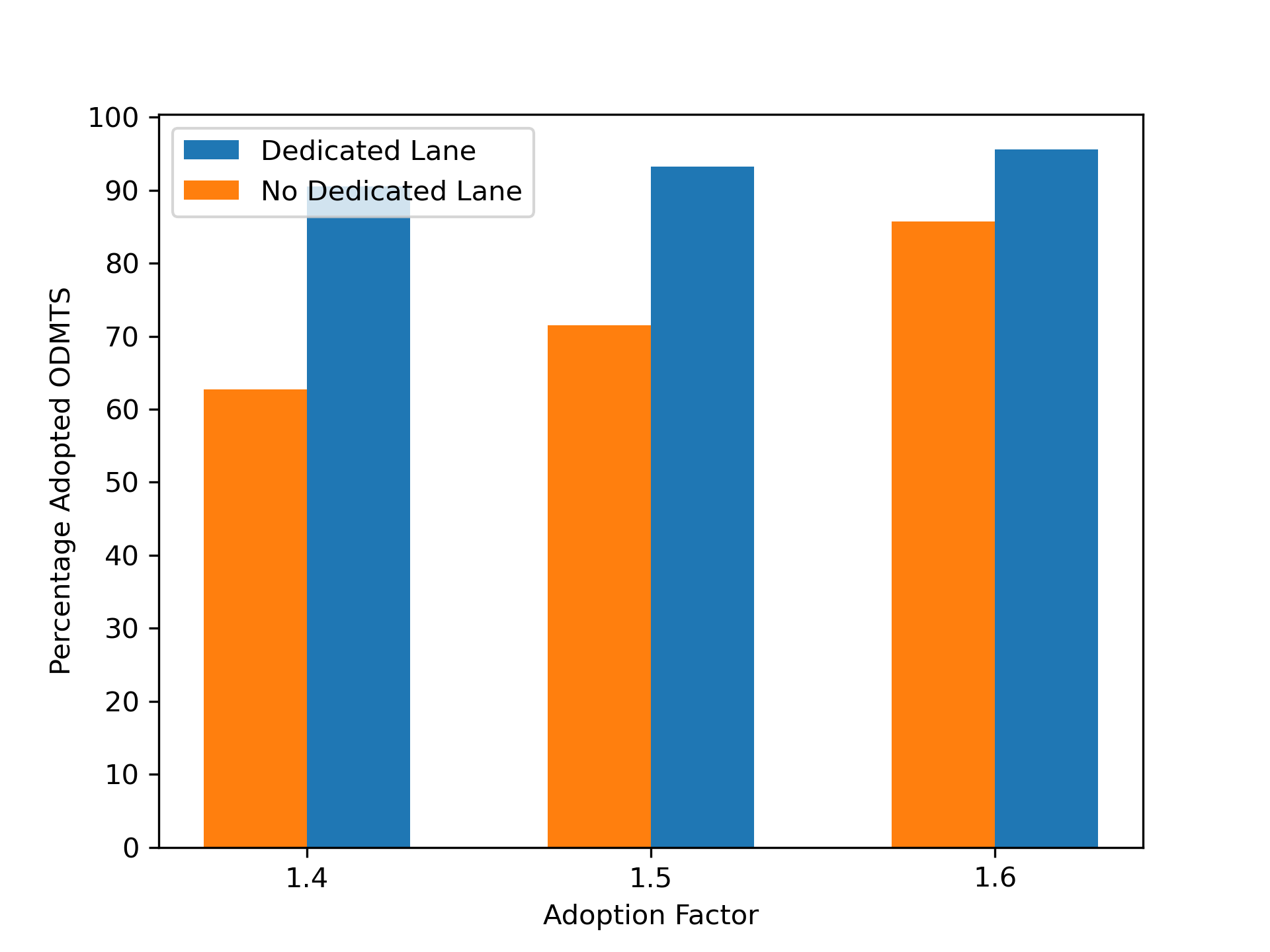}
		\caption{Pessimistic Scenario}
		\label{fig:adoption-percent-pessimistic}
	\end{subfigure}
	\caption{ODMTS Adoption Percentages for Non-local Potential Riders}
	\label{fig:adoption-percent}
\end{figure}

This section analyzes how changing the adoption factor affects
adoption by potential riders.  Recall that the adoption factor $\rho$
is used in the choice model to decide whether travelers adopt the
ODMTS: Potential riders adopt the system if the ODMTS travel time is
at most $\rho$ times the travel time by car.
Figure~\ref{fig:adoption-percent} shows the different adoption rates
for non-local potential riders when the network is designed for different scenarios and different
adoption factors ranging from $\rho=1.4$ (less willing to adopt) to
$\rho=1.6$ (more willing to adopt).  Even without DBLs, adoption rates
are typically high (above 50\%) across scenarios and adoption factors.
However, this is not the case for the expected scenario with adoption
factor $\rho=1.4$: Relatively low congestion favors driving by car.
This, combined with less willingness to adopt the ODMTS, results in an
adoption rate closer to 20\%.  Introducing DBLs in this situation has
a major effect on adoption, which climbs to above 60\%.  This
demonstrates the effectiveness of DBLs in attracting non-local riders.

\begin{figure}[t]
	\centering
	\includegraphics[width=.7\linewidth]{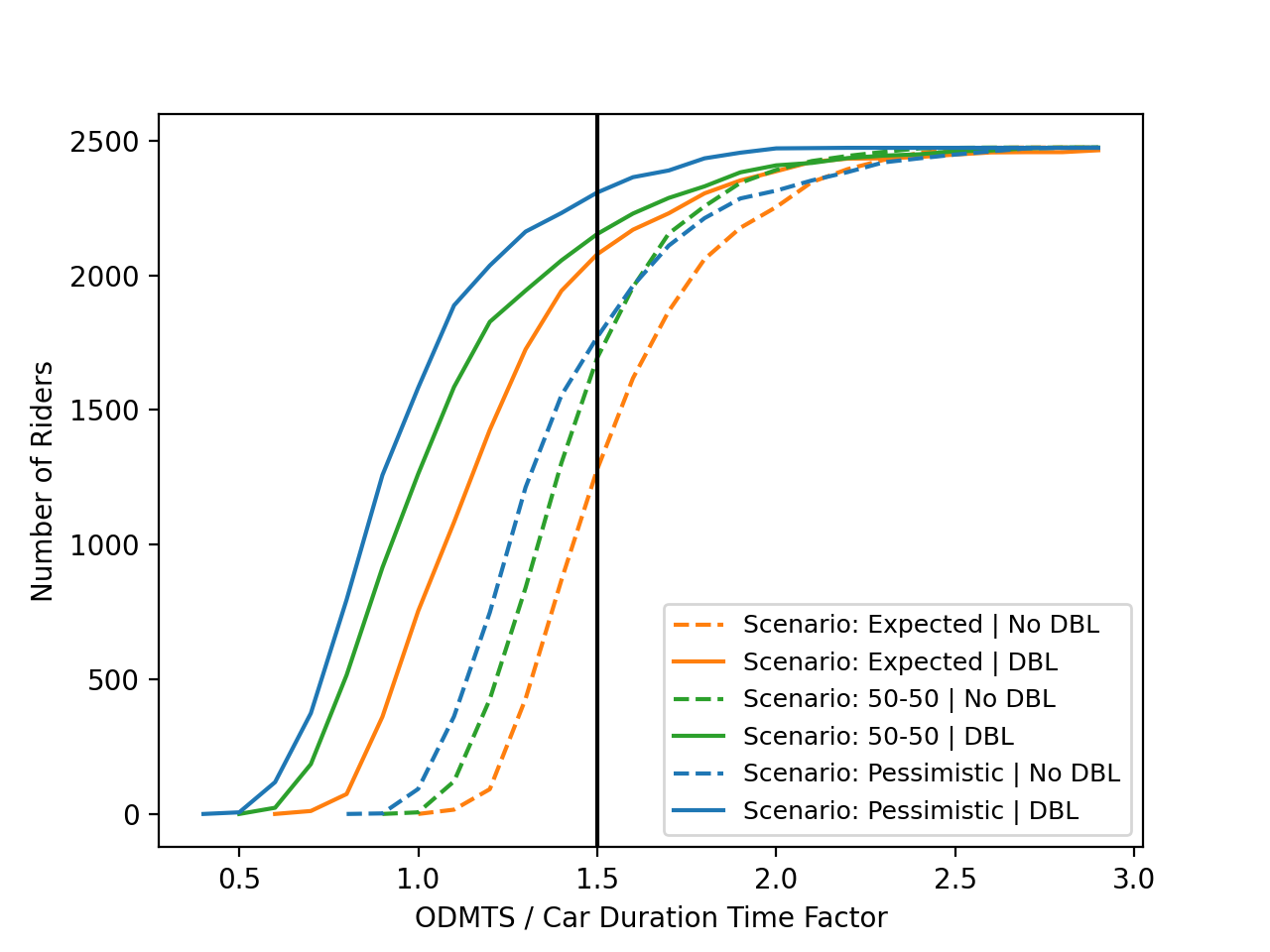}
	\caption{Ratio of ODMTS to Car Travel Times for Non-Local Potential Riders}
	\label{fig:cum_dist_odmts_vs_car}
\end{figure}

Figure~\ref{fig:cum_dist_odmts_vs_car} explains why DBLs are so
effective in increasing adoption.
It considers the networks designed for $\rho=1.5$ and shows the ratio of ODMTS to car
travel times for all non-local potential riders under the different
scenarios.
In other words, it shows the actual ratio that the potential riders observe, where riders adopt if the ratio is at most $\rho = 1.5$.
These curves are very steep.  For example, in the expected
scenario without DBLs ridership quickly increases around the $1.5$ cut-off. Including DBLs makes the
ODMTS more competitive and moves the curves to the left, resulting in
a sharp increase in adoption, especially when adoption was previously
low.

\section{Conclusion}
\label{sec:discussion}

This paper studied the effect of congestion on ODMTS and to what
extent Dedicated Bus Lanes (DBLs) can improve rider adoption.  To
enable this study, the paper introduced a method to create realistic
congestion scenarios without querying all origin-destination pairs by
combining data from multiple sources.  Designing the ODMTS was modeled
as a bilevel optimization problem in which potential riders choose
whether to adopt the system based on the quality of the trip offered.
These methods were then applied to perform a comprehensive case study
in the Metro Atlanta Area, based on real data by Gwinnett County
Transit and the Atlanta Regional Commission.

The case study revealed that, even with significant traffic, ODMTS is
able to provide accessible transit, experiences high adoption, and
almost breaks even in terms of cost.  That being said, Interstate 85
from Gwinnett to Atlanta is so congested that the ODMTS avoids using
this corridor at all.  Instead passengers are transfered to take a
detour through the rail system and experience long travel times.  The
congestion is not isolated to the highways, and slower shuttle trips
are a major source of increasing costs as congestion increases.

The paper then studied the case where DBLs were added to the network
and the ODMTS was reoptimized to benefit from the congestion-free
lanes.  The new networks feature buses on I-85 that are used
extensively by riders.  Non-local travel sees huge improvement in
terms of travel time and ODMTS adoption.  In fact, DBLs are so
effective that the travel times under the pessimistic scenario are
better than the expected scenario without DBLs.  These improvements do
not come at the expense of local travel, which sees a service similar
to the baseline.  While system costs still increase with congestion,
DBLs do not contribute negatively to the operating cost.

This case study shows that DBLs may significantly increase ODMTS
adoption.  Future work should compare the benefit of DBLs in ODMTS
across different cities, since Atlanta is a particularly congested
city relying heavily on its highway system. Within various cities, studies may also involve scenarios with different bus transit frequencies and different sets of existing riders. Another direction for
future research may be to improve the modeling of bus transit wait
times to more accurately assess travel times.

\section*{Acknowledgments}
This research is partly supported by NSF LEAP-HI proposal NSF-1854684.

\bibliographystyle{trc}
\bibliography{references}

\end{document}